\documentclass[a4paper, 12pt]{article}
\usepackage{amsmath,amsfonts,amsthm,amssymb,here,epsf}
\usepackage[latin1]{inputenc}
\usepackage[T1]{fontenc}
\usepackage{ae,aecompl}
\begin{document}

\newtheorem{theo}{Theorem}[section]
\newtheorem{prop}[theo]{Proposition}
\newtheorem{lemme}[theo]{Lemma}
\newtheorem{cor}[theo]{Corollary}
\newtheorem{defi}[theo]{Definition}
\newtheorem{assu}[theo]{Assumption}
\newtheorem{nontheo}[theo]{Conjectured theorem}


\newcommand{\beq}{\begin{eqnarray}}
\newcommand{\enq}{\end{eqnarray}}
\newcommand{\be}{\begin{eqnarray*}}
\newcommand{\en}{\end{eqnarray*}}
\newcommand{\T}{\Bbb T}
\newcommand{\N}{\Bbb N}
\newcommand{\Td}{\Bbb T^d}
\newcommand{\R}{\Bbb R}
\newcommand{\Rd}{\Bbb R^d}
\newcommand{\Zd}{\Bbb Z^d}
\newcommand{\Linf}{L^{\infty}}
\newcommand{\dt}{\partial_t}
\newcommand{\Dt}{\frac{d}{dt}}
\newcommand{\demi}{\frac{1}{2}}
\newcommand{\ep}{^{\epsilon}}
\newcommand{\epu}{_{\epsilon}}
\newcommand{\ds}{\displaystyle}
\newcommand{\bx}{{\mathbf X}}
\newcommand{\bsx}{{\mathbf x}}
\newcommand{\bg}{{\mathbf g}}
\newcommand{\bv}{{\mathbf v}}
\let\cal=\mathcal


\bibliographystyle{plain}


\title{A fully non-linear version of the incompressible Euler equations: the semi-geostrophic system}

\author{G. Loeper}

\maketitle


\begin{abstract}
This work gathers new results concerning the semi-geostrophic equations: existence and stability of measure valued solutions, existence and uniqueness of solutions under certain continuity conditions for the density, convergence to the incompressible Euler equations. Meanwhile, a general technique  to prove uniqueness of sufficiently smooth solutions to non-linearly coupled system
is introduced, using optimal transportation.
\end{abstract}

\tableofcontents

\section{Introduction}
The semi-geostrophic equations are an approximation to the Euler equations of fluid mechanics, used in meteorology to describe atmospheric flows, in particular they are believed (see \cite{CP}) to be an efficient model to describe frontogenesis. 
Different versions (incompressible \cite{BB2}, shallow water \cite{CuG} , compressible \cite{CuMa}) of this model have been studied, and we will
focus here on the incompressible 2-d and 3-d version.
The 3-d model describes the behavior of an incompressible fluid in a domain $\Omega\subset \R^3$. To the evolution in $\Omega$ is associated a motion in a 'dual' space, described by the following non-linear transport equation:
\be
&& \dt \rho + \nabla\cdot (\rho \bv)=0,\\
&& \bv=(\nabla\Psi(x)-x)^{\perp},\\
&& \det D^2\Psi = \rho,\\
&&\rho(t=0)=\rho^0.
\en
Here $\rho^0$ is a probability measure on $\R^3$, and for every $\bv=(v_1,v_2,v_3)\in \R^3$, $\bv^{\perp}$ stands for $(-v_2, v_1, 0)$. 
The velocity field is given at each time by solving a Monge-Amp\`ere equation 
in the sense of the polar factorization of maps (see \cite{Br1}), 
{\it i.e.} in the sense that $\Psi$ is convex from $\R^3$ to $\R$ and satisfies $\nabla\Psi_{\#}\rho= \chi_{\Omega}{\cal L}^3$, where ${\cal L}^3$ is the Lebesgue measure of $\R^3$, and $\chi_{\Omega}$ is the indicator function of $\Omega$. For compatibility $\Omega$ has Lebesgue-measure one.
This model arises as an approximation to the primitive equations of meteorology,
and we shall give a brief idea of the derivation of the model, although
the reader interested in more details should refer to \cite{CP}. 

In this work we will deal with various questions related to the semi-geostrophic (hereafter $SG$) system: existence and stability of measure-valued solutions, existence and uniqueness of smooth solutions, and finally convergence towards the incompressible Euler equations in 2-d. The results are stated in more details in section 1.

\subsection{Derivation of the semi-geostrophic equations}
We now give for sake of completeness  a brief and simplified idea of the derivation of the system, inspired from \cite{BB2}, and more complete arguments can be found in \cite{CP}.

\subsubsection*{Lagrangian formulation}
We start from the 3-d incompressible Euler equations with constant Coriolis parameter $f$ in a domain $\Omega$.
\be
&&\frac{D \bv}{Dt} + f \bv^{\perp} = \frac{1}{\rho}\nabla p - \nabla \varphi, \\
&& \nabla \cdot \bv =0, \ \frac{D \rho}{Dt}=0,\\
&& \bv\cdot \partial\Omega=0,
\en
where $\ds\frac{D \cdot}{Dt}$ stands for $\partial_t + \bv\cdot \nabla$, and we still use $\bv^{\perp}=(-v_2,v_1,0)$.
The term $\nabla\varphi$ denotes the gravitational effects (here we will take $\varphi=g x_3$ with constant $g$), and
the term $f\bv^{\perp}$ is the Coriolis force due to rotation of the Earth.
For large scale atmospheric flows, the Coriolis force $f\bv^{\perp}$  dominates the advection term $\ds\frac{D \bv}{Dt}$, and renders the flow mostly bi-dimensional.
We use the hydrostatic approximation: $\partial_{x_3}p=-\rho g$ and restrict ourselves to the case $\rho\equiv 1$. 

Keeping only the leading order terms leads to the geostrophic balance
\be
\bv_g=-f^{-1}\nabla^{\perp} p,
\en
that defines $\bv_g$, the geostrophic wind.
Decomposing $\bv=\bv_g+ \bv_{ag}$ where the second component is the ageostrophic wind, supposedly small departures from the geostrophic balance, the 
semi-geostrophic system reads:
\be
&&\frac{D \bv_g}{Dt} + f \bv^{\perp} = \nabla_H p,\\
&& \nabla \cdot \bv =0,
\en
where $\nabla_H=(\partial_{x_1}, \partial_{x_2},0)$.
Note however that the advection operator  $\partial_t + \bv\cdot \nabla$ still uses the full velocity $\bv$.
Introducing the potential 
\be
\Phi=\demi |x_H|^2 + f^{-2} p,
\en
with $x_H=(x_1,x_2,0)$, 
we obtain  the following
\be
\frac{D }{Dt}\nabla\Phi(t,x)=f (x-\nabla\Phi(t,x))^{\perp}.
\en 
We introduce the lagrangian map $\bg:\Omega\times \R^+ \mapsto \Omega$ giving the position at time $t$ of the parcel located at $x_0$ at time 0.
The previous equation means that, if for fixed $x\in\Omega$ we consider the trajectory in the 'dual' space, defined by
$X(t,x)=\nabla\Phi(t,\bg(t,x))$, we have 
\be
\dt X(t,x)=f\left(\bg(t,x)-X(t,x)\right)^{\perp}.
\en

By rescaling the time, we can set $f=1$. As stated the system looks under-determined: indeed $\Phi$ is unknown;
however we have the condition $X(t,x)=\nabla \Phi (t,\bg(t,x))$.
Moreover, the dynamic in the $x$ space being incompressible and contained in $\Omega$, the map $\bg(t,\cdot)$ must be measure preserving in $\Omega$ for each $t$, {\it i.e.} 
\be
{\cal L}^3(\bg(t)^{-1}(B))={\cal L}^3(B)
\en
for each $B\subset \Omega$ measurable (where ${\cal L}^3$  denotes the Lebesgue measure of $\R^3$). We shall hereafter denote by $G(\Omega)$ the set of all such measure preserving maps.
Then Cullen's stability criteria asserts that the potential $\Phi$ should be convex 
for the system to be stable to small displacements of particles in the $x$ space. Hence, for each $t$, $\Phi$ must be a convex function such that
\be
X(t,\cdot)=\nabla \Phi(t,\bg(t,\cdot)),
\en 
with $\bg(t,\cdot)\in G(\Omega)$. In the next paragraph we shall see that, under very mild assumptions on $X$, this decomposition, called polar factorization, can only happen
for a unique choice of $\bg$ and $\nabla \Phi$. Now if $\Phi^*$ is the Legendre transform of $\Phi$, 
\be
\Phi^*(y) = \sup_{x\in \Omega} x\cdot y - \Phi(x),
\en
then $\nabla \Phi$ and $\nabla \Phi^*$ are inverse maps of each other, and the semi-geostrophic system then reads
\be
&&\frac{D X}{Dt}= \left(\nabla \Phi^*(X(t))- X(t) \right)^{\perp},\\
&&\nabla \Phi^*(t)\circ X(t) \in G(\Omega).
\en

In the next paragraph, we expose the results concerning the existence and uniqueness of the gradients $\nabla \Phi,\nabla \Phi^*$.

\subsection{Polar factorization of vector valued maps}
The polar factorization of maps has been discovered by Brenier in \cite{Br1}. 
It has later been extended to the case of general Riemannian manifolds by McCann in \cite{Mc2}. 
\subsubsection*{The Euclidean case}
Let $\Omega$ be a fixed bounded domain of $\R^d$ of Lebesgue measure 1 and satisfying the 
condition ${\cal L}^d(\partial \Omega)=0$.
We consider a mapping $X\in L^2(\Omega;\Rd)$.
We will also consider the push-forward of the Lebesgue measure of $\Omega$  by $X$, that we will denote by $X_{\#} \chi_{\Omega}{\cal L}^d = d\rho$ (or, in short, $X_{\#} dx$) and which is defined by
\be
\forall f \in C^0_b(\Rd), \int_{\Rd} f(x) \ d\rho(x) = \int_{\Omega} f(X(x)) \ dx.
\en
Let  ${\cal P}$ be the set of  probability measures $\R^d$, and ${\cal P}^2_a$ the subset of ${\cal P}$
where the subscript $a$ means absolutely continuous with respect to the Lebesgue measure
 (or equivalently that have a density in $L^1(\Rd)$), and 
the superscript $2$ means 
with finite second moment. ({\it i.e.} such that
\be
\int_{\Rd} |x|^2 d \rho(x) < +\infty.)
\en
Note that for $X\in L^2(\Omega,\Rd)$, the measure $\rho=X{_{\#}} dx$ has necessarily finite second moment, and thus belongs to ${\cal P}^2$.
\begin{theo}[Brenier, \cite{Br1}]\label{5defpsi}
Let $\Omega$ be as above, $X\in L^2(\Omega;\Rd)$ and $\rho=X_{\#}dx$.
\begin{enumerate}
\item There exists a unique up to a constant convex function, that will be denoted $\Phi[\rho]$, such that:
\be
\forall f \in C^0_b(\Rd),\ \int_\Omega f(\nabla\Phi[\rho](x)) \ dx =\int_{\Rd} f(x)d\rho(x).
\en
\item Let $\Psi[\rho]$ be the Legendre transform of $\Phi[\rho]$,
if $\rho\in {\cal P}^2_a$, $\Psi[\rho]$
 is the unique up to a constant convex function satisfying 
\be
\forall f \in C^0_b(\Omega),\ \int_{\Rd} f(\nabla\Psi[\rho](x)) \ d\rho(x) =\int_\Omega f(x) dx.
\en
\item If $\rho\in {\cal P}^2_a$, 
$X$ admits the following unique polar factorization:
\be
X=\nabla\Phi[\rho]\circ g,
\en
with $\Phi[\rho]$ convex, $g$ measure preserving in $\Omega$.
\end{enumerate}
\end{theo}

{\it Remark:} $\Psi[\rho], \Phi[\rho]$ depend only on $\rho$, and are solutions (in some weak sense) 
respectively in $\Rd$ and $\Omega$, of the Monge-Amp\`ere equations 
\be
&&\det D^2\Psi=\rho,\\
&&\rho(\nabla\Phi) \ \det D^2\Phi =1.
\en
When $\Psi$ and $\Phi$ are not in $C^2_{loc}$ these equations can be understood in the viscosity (or Alexandrov)
sense or in the sense of Theorem \ref{5defpsi}, which is strictly weaker.
For the regularity of those solutions and the consistency of the different weak formulations
the reader can refer to \cite{Ca3}.

\subsubsection*{The periodic case}
The polar factorization theorem has been extended to Riemannian manifolds in \cite{Mc2} (see also \cite{Co} for the case of the flat torus). In this case, we consider a mapping $X: \Rd \mapsto \Rd$ such that for all $\vec p \in \Zd, X(\cdot + \vec p) = X+\vec p$. Then $\rho=X_{\#}dx$ is a probability measure on $\Td$. 
We define $\Psi[\rho], \Phi[\rho]$ through the following:
\begin{theo}\label{5defpsiper}
Let $X:\Rd \to \Rd$ be as above, with $\rho=X_{\#}dx$. 
\begin{enumerate}
\item Up to an additive
 constant there exists a unique
convex function $\Phi[\rho]$ such that
$\Phi[\rho](x)-x^2/2$ is \mbox{$\Zd$-periodic} (and thus $\nabla\Phi[\rho](x)-x$ is $\Zd$ periodic), and 
\be
\forall f \in C^0(\Td),\ \int_{\Td} f(\nabla\Phi[\rho](x)) \ dx =\int_{\Td} f(x)  \ d\rho(x).
\en
\item Let $\Psi[\rho]$ be the Legendre transform of $\Phi[\rho]$. If $\rho$ is Lebesgue integrable, $\Psi[\rho]$ is the unique up to a constant convex function satisfying
\\
$\Psi[\rho](x)-x^2/2$ is \mbox{$\Zd$-periodic} (and thus $\nabla\Psi[\rho](x)-x$ is $\Zd$ periodic), and 
\be
\forall f \in C^0(\Td),\ \int_{\Td} f(\nabla\Psi[\rho](x)) \ d\rho(x) =\int_{\Td} f(x) \ dx.
\en
\item If $\rho$ is Lebesgue integrable, $X$ admits the following unique polar factorization:
\be
X=\nabla\Phi[\rho] \circ g
\en
with $g$ measure preserving from $\Td$ into itself, and $\Phi[\rho]$ convex, $\Phi[\rho]-|x|^2/2$ periodic.
\end{enumerate}
\end{theo}

{\it Remark 1:} From the periodicity of $\nabla\Phi[\rho](x)-x, \nabla\Psi[\rho](x)-x$, for every $f$ \mbox{$\Zd$-periodic}, $f(\nabla\Psi[\rho]), f(\nabla\Phi[\rho])$ are well 
defined on $\Rd/\Zd$.

{\it Remark 2:} Both in the periodic and non periodic case, the definitions of $\Psi[\rho]$ and $\Phi[\rho]$  make sense if $\rho$ is absolutely continuous with respect to the Lebesgue measure. If not, the definition and uniqueness of $\Phi[\rho]$ is still valid,
as well as the property $\nabla\Phi[\rho]_{\#}\rho=\chi_{\Omega}{\cal L}^d$.
The definition of $\Psi[\rho]$ as the Legendre transform of $\Phi[\rho]$ is still valid also, but then the expression $\ds\int f(\nabla\Psi[\rho](x)) \ d\rho(x)$ does not necessarily make sense  since $\nabla\Psi$ is not necessarily continuous. Moreover the polar factorization does not hold any more.

{\it Remark 3:} We have (see \cite{Co}) the unconditional bound $$\|\nabla\Psi[\rho](x)-x\|_{\Linf(\Td)}\leq \sqrt{d}/2$$
that will be useful later on. 

\subsection{Lagrangian formulation of the $SG$ system}
From Theorems \ref{5defpsi}, \ref{5defpsiper} 
the Lagrangian formulation of the semi-geostrophic equation then becomes
\beq\label{5lagrangian}
&&\frac{D X}{Dt}= \left[\nabla\Psi(X) -X \right]^{\perp},\\
&&\Psi=\Psi[\rho], \  \  \rho=X{_{\#}}dx.\label{5lagrangian2}
\enq

\subsection{Eulerian formulation in dual variables}
In both cases (periodic and non periodic) we thus investigate the following system that will be referred to as $SG$: we look for a time dependent probability measure $t\to \rho(t,\cdot)$ satisfying 
\begin{eqnarray}
&&\partial_t \rho  +\nabla\cdot (\rho \bv)= 0\label{5sg1},\\
&&\bv(t,x)=\left(\nabla\Psi[\rho(t)](x)-x\right)^{\perp}\label{5sg2},\\
&&\rho(t=0)=\rho^0\label{5sg3}.
\end{eqnarray}
Weak solutions (which are defined below) of this system with $L^p$ initial data for $p\geq 1$ have been found, see
\cite{BB2}, \cite{CuG}, \cite{LN}.

\subsection{Results}
In this work we deal with various mathematical problems related to this system: we first extend the notion of weak solutions
that had been shown to exist for $\rho\in \Linf(\R_+, L^q(\R^3))$, $q>1$ (\cite{BB2}, \cite{CuG}), and then for  $\rho\in \Linf(\R_+, L^1(\R^3))$ (\cite{LN}), to the more general case of bounded measures. 
The question of existence of measure-valued solutions was raised and left unanswered in those papers, and we show here existence of global solutions to the Cauchy problem with initial data a bounded compactly supported measure, and show the weak stability/compactness of these {\it weak measure} solutions.

Then we show existence of continuous solutions, more precisely, we show local existence of solutions with Dini-continuous (see (\ref{5modulus})) density. 
For this solutions, the velocity field is then $C^1$ and the Lagrangian system (\ref{5lagrangian},\ref{5lagrangian2}) is defined  everywhere.

We also show uniqueness in the class of H\"older continuous solutions (a sub-class of Dini continuous solutions).
This proof uses in an original way the optimal transportation of measures by convex gradients and its regularity properties,
and can be adapted to give a new proof of uniqueness for solutions of the 2-d Euler equation with bounded vorticity, but also for a broad class of non-linearly coupled system. The typical application is a density evolving through a transport equation where the velocity field depends on the gradient of a potential, the potential solving an elliptic equation with right hand side the density. Well known examples of such cases are the Vlasov-Poisson and Euler-Poisson systems.

Finally, in the 2-d case, we study the convergence of the system to the Euler incompressible equations; this convergence is expected for $\rho$ close to 1, since formally expanding $\Psi=x^2/2 + \epsilon\psi$, and linearizing the determinant around the identity matrix,  we get
$$\det D^2\Psi=1 + \epsilon \Delta\psi + O(\epsilon^2),$$ and the Monge-Amp\`ere equation turns into the Poisson equation  
$$\Delta\psi=\frac{\rho-1}{\epsilon}=:\mu.$$
After a proper time scaling, $\mu$ satisfies 
\be
&&\dt\mu + \nabla\cdot(\mu \nabla^{\perp}\psi),\\
&&\Delta \psi =\mu,
\en 
that we recognize as the vorticity formulation of the 2-d Euler incompressible equation.
The study of this 'quasi-neutral' limit is done by two different ways: One uses a modulated energy method similar as the one used in \cite{Br2} and \cite{BL} and is valid for weak solutions. The other uses a more classical expansion of the solution, and regularity estimates, and is similar to the method used in \cite{L1}. The second method yields also a time of existence for the smooth solution that goes to infinity, as the scaling parameter $\epsilon$ goes to 0.
From a physical point of view, this asymptotic study may be seen as a justification of the consistency of the semi-geostrophic approximation.

\section{Measure valued solutions}

\subsection{A new definition of weak solutions}
We have first the following classical weak formulation of equation (\ref{5sg1}):
\\ 
$\rho \in C(\R+, L^1(\R^3)-w)$ is said to be a weak solution of $SG$ if  
\be
&&\forall T>0,  \ \forall \varphi \in C^{\infty}_c([0,T]\times \Bbb R^2),\\ 
&&\int \partial_t \varphi \,\rho  +\nabla\varphi \cdot (\nabla\Psi[\rho]-x)^{\perp}\,\rho\;dtdx\\
&=& \int \varphi(T,x)\rho(T,x)dx -\int \varphi(0,x)\rho(0,x)dx,
\en
where for all $t$, $\Psi[\rho]$ is as in Theorem \ref{5defpsi}.
The problematic part in the case of measure valued solutions is to give sense to the product
$\rho \nabla\Psi[\rho]$ since at the point where $\rho$ is singular $\nabla\Psi[\rho]$ is unlikely to be continuous.
Therefore we use the Theorem \ref{5defpsi} to write for any $\rho \in {\cal P}^2_a(\R^3)$
\be
\forall \varphi \in C^{\infty}_c(\R^3), \int_{\R^3} \rho \nabla\Psi[\rho]^\perp \cdot \nabla\varphi=\int_\Omega x^\perp\cdot \nabla\varphi(\nabla\Phi[\rho])
\en
(the integrals would be performed over $\T^3$ in the periodic case). 
The property $\nabla\Phi[\rho]_{\#} \chi_{\Omega}{\cal L}^3 = \rho$  is still valid when $\rho$ is only a measure with finite second moment (see Remark 2 after Theorem \ref{5defpsiper}). Therefore, the formulation on the right hand side extends unambiguously to the case where $\rho \notin L^1(\R^2)$. 
\subsubsection*{Geometric interpretation}
This weak formulation allows has a natural geometric interpretation:
at a point where $\Psi[\rho]$ is not differentiable, and 
thus where $\partial\Psi[\rho]$ is not reduced to a single point, 
$\nabla\Psi[\rho]$ should be replaced by $\bar{\partial}\Psi[\rho]$
the center of mass of the (convex) set $\partial\Psi[\rho]$.
\\
This motivates the following definition of  weak measure solutions
\begin{defi}\label{5weak}
Let, for all $t\in [0,T]$, $\rho(t)$ be a probability measure of $\R^3$. It is said to be a {\bf weak measure solution} to $SG$ with initial data $\rho^0$ if 

\bigskip
\indent 1-
The time dependent probability measure $\rho$  belongs to $C([0,T], {\cal P}-w*)$,

\bigskip
\indent 2- 
there exists $t\to R(t)$  non-decreasing such that
for all $t\in [0,T]$, $\rho(t,\cdot)$ is supported in $B(0,R(t))$,
 
\bigskip
\indent 3- for all $T>0$  and for all $\varphi \in C^{\infty}_c([0,T]\times \Bbb R^3)$ we have
\begin{eqnarray}\label{5eqweak}
&&  \  \  \ \int_{[0,T]\times \R^3} \partial_t \varphi(t,x)  \ d\rho(dt,x)  \\
&+&\int_{[0,T]\times \Omega}\nabla\varphi(t,\nabla\Phi[\rho(t)](x)) \cdot x^{\perp} \ dtdx-\int_{[0,T]\times \R^3}\nabla\varphi(t,x) \cdot  x^{\perp} \ d\rho(dt,x) \nonumber  \\
&=&\int \varphi(T,x)d \rho(T,x) \ dx - \int \varphi(0,x)d \rho^0(x) \ dx\nonumber. 
\end{eqnarray}
\end{defi}
This definition is consistent with the classical definition of weak solutions if for all $t$, $\rho(t,\cdot)$ is absolutely continuous with respect to the Lebesgue measure.

\subsection{Result} Here we prove the following
\begin{theo}\label{5main}
\begin{enumerate}
\item Let $\rho^0$ be a probability measure compactly supported.  There exists
a global weak measure solution to the system $SG$ with initial data $\rho^0$ in the sense of Definition \ref{5weak}.
\item For any $T>0$, if $(\rho_n)_{n\in\N}$ is a sequence of weak measure solutions on $[0,T]$ to $SG$ with initial data $(\rho^0_n)_{n\in\N}$, supported in $B_R$ for some $R>0$ independent of $n$, 
 the sequence $(\rho_n)_{n\in \N}$ is precompact in $C([0,T], {\cal P}-w*)$ and every converging subsequence converges to a weak measure solution of $SG$.
\end{enumerate}
\end{theo}

\subsection*{Proof of Theorem \ref{5main}}

We first show the weak stability of the formulation of Definition (\ref{5weak}), and the compactness of weak measure solutions. We then use this result to obtain global existence of solutions to the Cauchy problem with initial data a bounded measure. 

\subsubsection*{Weak stability of solutions}

We consider a sequence $(\rho_n)_{n\in \N}$ of solutions of $SG$ in the sense of Definition \ref{5weak}. The sequence is uniformly compactly supported at time 0. 
We first show that there exists a non-decreasing function $R(t)$ such that
$\rho_n(t)$ is supported in $B(R(t))$ for all $t,n$:
\begin{lemme}
Let $\rho\in C([0,T], {\cal P}(\R^3)-w*)$ satisfy (\ref{5eqweak}), let \mbox{$\rho^0=\rho(t=0)$} be supported in $B(0,R^0)$, then $\rho(t)$ is supported in 
$B(0, R^0 + Ct) $, \mbox{$C=\sup_{y\in \Omega}\{|y|\}$}.
\end{lemme} 
{\it Proof.} Consider any  function $\xi\epu(t,r)\in C^{\infty}(\R)$ such that 
\be
&&\xi\epu(0,r)\equiv 1 \ {\rm if } \  -\infty <  r  \leq R^0,    \\
&&\xi\epu(0,r)\equiv 0 \ {\rm if } \  r \geq R^0+\epsilon,  \\
&&\xi\epu(t,r)=\xi\epu(r-Ct),
\en 
with $\xi(0,\cdot)$ non increasing. 
Then compute 
\be
&&\Dt \int \xi\epu(t,|x|) \  d\rho(t,x)\\
&=& -\int \partial_r\xi\epu(t,|x|)C \  d\rho(t,x)
+ \int_{\Omega}  \partial_r\xi\epu(t,|\nabla\Phi[\rho(t)]|)\frac{\nabla\Phi[\rho(t)]}{|\nabla\Phi[\rho(t)]|}\cdot x^{\perp} \ dx\\
&\geq& \int_{\Omega}\partial_r\xi\epu(t,|\nabla\Phi[\rho(t)]|)(-C+ |x|) \ dx\\
&\geq& 0
\en
since, by definition, for $x\in\Omega$, $|x|\leq C$ and $\xi$ is non increasing with respect to $r$.
Note also that we have used  
$\ds\int \nabla_x[\xi(t,|x|)] \cdot x^{\perp} d\rho(t,x) \ dx \equiv 0$.
We know on the other hand that 
\be
&&\int_{\R^3} \xi\epu(0,|x|) d\rho(0,x)=1,\\
&&\int_{\R^3} \xi\epu(t,|x|) d\rho(t,x)\leq 1,
\en
therefore we conclude that $\ds\int_{\R^3} \xi\epu(|x|, t) d\rho(t,x)\equiv 1$,
which concludes the lemma by letting $\epsilon$ go to 0.

$\hfill\Box$

From this lemma, we have: $$\left|-\int_{[0,T]\times \R^3}\nabla\varphi(t,x) \cdot  x^{\perp} \ d\rho_n(dt,x) \nonumber +  \int_{[0,T]\times \Omega}\nabla\varphi(t,\nabla\Phi[\rho_n(t)](x)) \cdot x^{\perp} \ dtdx\right|$$
$$\leq C(T) \|\varphi\| _{L^1([0,T], C^1(B_{R(T)})}.$$
Thus from Definition \ref{5weak} equation (\ref{5eqweak})  we  know that for any time $t\geq 0$, $\dt\rho_n(t,\cdot)$ is bounded in the dual of $L^1([0,T],C^1(\R^3))$ and thus in the dual of $L^1([0,T],W^{2,p}(\R^3))$ for $p>3$ by Sobolev embeddings. Thus for some $p'>1$ we have
\be
\dt\rho_n \in \Linf([0,T], W^{-2,p'}(\R^2)).
\en

With the two above results, and using classical arguments of functional analysis (see \cite{Li}), we can obtain the following lemma:
\begin{lemme}
Let the sequence $(\rho_n)_{n\in\N}$ be as above, there exists $\rho \in C([0,T], {\cal P}-w*)$ and a subsequence $(\rho_{n_k})_{k\in\N}$, such that for all $t\in [0,T]$,  $\rho_{n_k}(t)$ converges to $\rho(t)$ in the weak-$*$ topology of measures.
\end{lemme}
With this lemma, we need to show that for all $\varphi \in C^{\infty}_c([0,T]\times \R^3)$ we have $\nabla\varphi(t, \nabla\Phi[\rho_n(t)])$ converging to $\nabla\varphi(t, \nabla\Phi[\rho(t)])$ whenever $\rho_n(t)$ converges weakly-$*$ to $\rho(t)$. This last step will be a consequence of the following stability theorem:
\begin{theo}[Brenier, \cite{Br1}]\label{5convergence}
Let $\Omega$ be as above. Let $(\rho_n)_{n\in\N}$ be a sequence of probability measures on $\Rd$, such that 
$\forall n,\;\int  (1+|x|^2)d\rho_n\leq C$,   let $\Phi_n=\Phi[\rho_n]$ and 
$\Psi_n=\Psi[\rho_n]$ be as in Theorem \ref{5defpsi}. 
If for any $f\in C^0(\Rd)$ such that $|f(x)|\leq C(1+|x|^2),$   $\int f\rho_n\rightarrow \int\rho f,$ then $\Phi_n\rightarrow \Phi[\rho]$
uniformly on each compact set of $\Omega$ and strongly in $W^{1,1}(\Omega;\Rd),$
and  $\Psi_n\rightarrow \Psi[\rho]$ uniformly on each compact set of $\Rd$ and strongly in $W^{1,1}_{loc}(\Rd)$. 
\end{theo} 
From this result, we obtain that the sequence
$\nabla\Phi[\rho_{n}]$ converges strongly in $L^1(\Omega)$ and almost everywhere (because of the convexity of $\Phi[\rho]$) to $\nabla\Phi[\rho]$. 
Thus $\nabla\varphi(t, \nabla\Phi[\rho_n])$ converges to  $\nabla\varphi(t, \nabla\Phi[\rho])$ in $L^1(\Omega)$ and one can pass to the limit in the formulation of Definition \ref{5weak}. This ends the proof of point 2 of Theorem \ref{5main}. 

\subsubsection*{Existence of solutions}
We show briefly the existence of a solution to the Cauchy problem in the sense of Definition \ref{5weak}. Indeed given $\rho^0$ the initial data for the problem that we want to solve, by smoothing $\rho^0$, we can take a sequence $\rho_n^0$ of initial data belonging to 
$L^1(\R^2)$, uniformly compactly supported and  converging \mbox{weakly-$*$} to $\rho^0$.
We know already from \cite{BB2}, \cite{CuG}, \cite{LN} that for every $\rho_n^0$,
one can build a global weak solution of (\ref{5sg1}, \ref{5sg2}, \ref{5sg3}), that will be uniformly
compactly supported on $[0,T]$ for all $T\geq 0$. This sequence will also
be solution in the sense of Definition \ref{5weak}. We then use the stability result, and conclude that, up to extraction of a subsequence, the sequence $\rho_n$ converges in $C([0,T], {\cal P}-w*)$ to a weak measure solution of $SG$  
with initial data $\rho^0$.
This achieves the proof of Theorem \ref{5main}. $\hfill\Box$

{\it Remark:} One can prove in fact the more general result, valid for non linear functionals:
\begin{prop}
Let $F\in C^0(\Omega\times \Rd)$, such that $|F(x,y)|\leq C(1+|y|^2)$, 
let $(\rho_n)_{n\in\N}$ be a bounded sequence of probability measures, Lebesgue integrable, with 
finite second moment. Let $\rho$ be a probability measure with finite second moment, such that 
for all $f \in C^0(\R^d)$ such that $|f(x)| \leq C(1+|x|^2)$,
$ \ds\int f d\rho_n \to \int f d\rho$. 
Then as $n$ goes to $\infty$, we have
\be
&&\int_{\Rd} F(\nabla\Psi[\rho_n](x),x) \ d\rho_n(x)
= \int_{\Omega} F(y,\nabla\Phi[\rho_n](y)) \ dy \\
&\to_{n}& \int_{\Omega} F(y,\nabla\Phi[\rho](y)) \ dy :=
\int_{\Rd} F(\bar{\partial}\Psi[\rho](x),x) \ d\rho(x).
\en
\end{prop}

\section{Continuous solutions}
What initial regularity is necessary in order to guarantee that the velocity fields remains Lipschitz, or that the flow remains continuous, at least for a short time ?
The celebrated Youdovich's Theorem for the Euler incompressible equation 
shows that when $d=2$, if the initial vorticity data is bounded in $\Linf$, 
the flow is H\"older continuous, with H\"older index decreasing to 0 as time goes to infinity.  This proof relies on the following regularity property of the Poisson equation: if $\Delta \phi$ is bounded in $\Linf$, then 
$\nabla\phi$ is Log-Lipschitz. This continuity is enough to define
a H\"older continuous flow for the vector field $\nabla\phi^\perp$. 
Such a result is not valid for the Monge-Amp\`ere equation. 
As far as we know, the optimal regularity result for Monge-Amp\`ere equations is the following:

\subsection{Regularity of solutions to Monge-Amp\`ere equation with Dini-continuous right hand side}
\begin{theo}[Wang, \cite{W}]\label{5wang}
Let $u$ be a strictly convex Alexandrov solution of 
\beq\label{5ma}
\det D^2 u = \rho
\enq
with $\rho$ strictly positive.
If  $w(r)$, the modulus of continuity of $\rho$,  satisfies 
\beq\label{5modulus}
\int_{0}^{1}\frac{w(r)}{r}dr < \infty,
\enq
 then $u$ is in $C^2_{loc}$.
\end{theo}
We will work here in the periodic case. In this case, $u$ the solution of (\ref{5ma}) will be $\Psi[\rho]$ of Theorem \ref{5defpsiper}.
The arguments of
 \cite{Ca0}, \cite{Ca3}, adapted to the periodic case, show that $\Psi[\rho]$ is indeed a strictly convex Alexandrov solution of solution of (\ref{5ma}).
Therefore we obtain the following corollary of Theorem \ref{5wang}:
\begin{cor}\label{5corwang} Let $\rho\in {\cal P}(\T^d)$  be such that 
\be
&&0<m\leq \rho \leq M,\\
&&\int_{0}^{1}\frac{w(r)}{r}dr = C < \infty. 
\en
where $m, M, C$ are positive constants.
Let $\Psi[\rho]$ be as in Theorem \ref{5defpsiper}. We have, for some constant ${\cal H}$ depending on $m,M,C$ 
\be
\|\Psi[\rho]\|_{C^2(\Td)} \leq {\cal H}.
\en
\end{cor}

\subsection{Result}
We will now prove the following:
\begin{theo}\label{5main3}
Let $\rho^0$ be a probability on $\T^3$, such that $\rho$ is  strictly positive and satisfies the continuity condition (\ref{5modulus}).
Then there exists $T>0$ and $C_1, C_2$ depending on $\rho^0$, such that on $[0,T]$
there exists a solution $\rho(t,x)$ of $SG$ that satisfies for all $t\in [0,T]$:  
\be
&&\int_{0}^{1}\frac{w(t,r)}{r}dr \leq C_1,\hspace{1cm} \|\Psi(t,\cdot)\|_{C^{2}(\T^3)} \leq C_2,
\en
where 
$w(t,r)$ is the modulus of continuity (in space) of $\rho(t,.)$.
\end{theo}

\subsection*{Proof of Theorem \ref{5main3}}
Let us first sketch the proof:
If $\Psi \in C^2$, then the flow $t\rightarrow X(t,x)$ generated by the velocity field $[\nabla\Psi(x)-x]^{\perp}$ is Lipschitz in space.
Since the flow is incompressible, we have
$\rho(t,x)=\rho^0(X^{-1}(t,x))$. 
\\
Now we use the following property: 
If two functions $f,g$ have modulus of continuity respectively $w_f, w_g$
then $g\circ f$ has modulus $w_g \circ w_f$. 
\\
Thus if $X^{-1}(t)$ is Lipschitz, we have  $w_{\rho^0\circ X^{-1}(t)}\leq w_{\rho^0}(L \ \cdot)$ with $L$ the Lipschitz constant of $X^{-1}(t)$ and condition (\ref{5modulus}) remains satisfied.

{\it Remark 1:} Note that H\"older continuous functions satisfy the condition (\ref{5modulus}).

{\it Remark 2:} Note also that we do  not need any integrability on $\nabla\rho$ and the solution of the Eulerian system still has to be understood in the distributional sense.

\subsubsection*{A fixed point argument}
Let us introduce the semi-norm 
\beq\label{5defnorme}
\|\mu\|_{\cal C} =   \int_{0}^{1}\frac{w_{\mu}(r)}{r}dr
\enq
defined on ${\cal P}(\T^3)$, where we recall that $w_{\mu}$ is the modulus of continuity of $\mu$. 
We denote ${\cal P_C}$ the set ${\cal P}$ equipped with this semi-norm, {\it i.e.}
\be
{\cal P_C}=\{ \mu \in {\cal P}(\T^3), \|\mu\|_{\cal C}<\infty \}.
\en
From now, we fix $\rho^0$ a probability density in ${\cal P_C}$,
satisfying $m\leq \rho^0 \leq M$, where $m$ and $M$ are strictly positive constants.  
Let $\mu$ be a time dependent probability density in $\Linf([0,T];{\cal P_C})$, such that $m\leq \mu(t)\leq M$ for all $t$,  
we consider  the solution $\rho$ of the initial value problem:
\beq
&&\dt\rho + (\nabla\Psi[\mu](x)-x)^\perp \cdot \nabla \rho=0\label{5lineaire},\\
&&\rho(t=0)=\rho^0\label{5lineaire-init}.
\enq
From Theorem \ref{5wang} and its corollary, the vector field $\bv[\mu]=(\nabla\Psi[\mu](x)-x)^{\perp}$ is $C^1$ uniformly in time, therefore there exists a unique solution to this equation, by Cauchy-Lipschitz Theorem.
This solution can be built by the method of characteristics as follows:
Consider the flow $X(t,x)$ of the vector field $\bv[\mu]$, then $\rho(t)$ is $\rho^0$ pushed forward by $X(t)$, {\it i.e.} $\rho(t)= \rho^0\circ X^{-1}(t)$.
 From the incompressibility of $\bv[\mu]$ the condition $m\leq\rho^0\leq M$ implies that for all $t\in [0,T]$, $m \leq \rho(t) \leq M$.

The initial data $\rho^0$ being fixed, the map $\mu\mapsto \rho$ will be denoted by ${\cal F}$.

The spatial derivative of $X$,  $D_xX$ satisfies 
\be
\dt D_x X= D_x \bv[\mu](X) D_x X,
\en
therefore we have
\be
\vert D_x X(t) \vert \leq \exp(t\sup_{s\in [0,t]}\vert D_x \bv[\mu](s) \vert),
\en
and the same bound holds for $X(t)^{-1}$. 
Since $w_{f\circ g} \leq  w_f \circ w_g$, and writing $C_t= \exp(t\sup_{s\in [0,t]}\vert D_x \bv[\mu] \vert)$, we obtain $w_{\rho(t)}(\cdot)\leq w_{\rho^0}(C_t \cdot)$, and 
\be
\int_0^1 \frac{w_{\rho(t)}(r)}{r}dr&\leq& \int_0^{C_t}\frac{w_{\rho^0}(r)}{r}dr\\
&\leq&\int_0^{1}\frac{w_{\rho^0}(r)}{r}dr + (M-m) (C_t-1),
\en 
(using that $\forall r, w_{\rho}(r) \leq M-m$).
Therefore, $$\|\rho(t)\|_{\cal C} \leq \|\rho^0\|_{\cal C} + (M-m)(C_t-1).$$
Now from Corollary \ref{5corwang}, and $m, M$ being fixed,  there exists a non-decreasing function ${\cal H}$ such that
\be
\|\bv[\mu]\|_{C^{1}} \leq {\cal H}(\|\mu\|_{\cal C}),
\en
and so $\ds C_t \leq \exp(t {\cal H}(\|\mu\|_{\Linf([0,t];{{\cal P_C}}})))$.
Hence we can chose $Q>1$, and then $T$ such that
\be
\|\rho^0\|_{\cal C} + (M-m)\left(\exp(T \  {\cal H}(Q\|\rho^0\|_{{\cal C}}))-1\right) = Q \|\rho^0\|_{\cal C}.
\en 
Note that for $Q>1$, we necessarily have $T>0$. 
Then  the map ${\cal F}:\mu \mapsto \rho$ goes now from
$$ {\cal A}= \left\{\mu, \ \|\mu\|_{\Linf([0,T]; {\cal P_C})} \leq Q \|\rho^0\|_{\cal C}, \  m\leq \mu \leq M \right\}$$
 into 
$$ {\cal B}= \left\{\rho, \ \|\rho(t)\|_{{\cal C}} \leq \|\rho^0\|_{\cal C} + (M-m)\left(\exp(t \  {\cal H}(Q\|\rho^0\|_{{\cal C}}))-1 \ \right), \forall t\in [0,T]  \right\},$$ and with our choice of $T=T(Q)$, we have ${\cal B} \subset {\cal A}$. 
Moreover from the unconditional bounds
\be
&&\rho \leq M,\\
&& \|\bv[\mu]\|_{\Linf([0,T]\times\T^3)}\leq \sqrt{3}/2,
\en
(see the remark after Theorem \ref{5defpsiper} for the second bound)
and using equation (\ref{5lineaire}),
we have also $\|\dt \rho\|_{\Linf([0,T]; W^{-1,\infty})} \leq K(M)$ whenever $\rho ={\cal F}(\mu)$.

Call $\tilde {\cal A}$ (resp. $\tilde{\cal B}$) the set 
${\cal A} \cap \{\rho, \|\dt \rho\|_{\Linf([0,T]; W^{-1,\infty})} \leq K(M) \}$,
(resp. ${\cal B} \cap \{\rho, \|\dt \rho\|_{\Linf([0,T]; W^{-1,\infty})} \leq K(M)\}$); we claim that 
\begin{itemize}
\item ${\cal F}(\tilde{\cal A})\subset \tilde{\cal B} \subset \tilde{\cal A}$,
\item $\tilde{\cal A}$ is convex and compact for the $C^0([0,T]\times \T^3)$ topology,
\item ${\cal F}$ is continuous for this topology,
\end{itemize}
so that we can apply the Schauder fixed point Theorem.
We only check the last point, the second being a classical result of functional analysis.
So let us consider a sequence $(\mu_n)_{n\in \N}$ converging to $\mu\in {\cal A}$, and the corresponding sequence \mbox{$(\rho_n={\cal F}(\mu_n))_{n\in\N}$.}
The sequence $\rho_n$ is pre-compact in $C^0([0,T]\times \T^3)$, from the previous point, and we see (with the stability Theorem \ref{5convergence}) that it converges to a solution $\rho$ of 
\be
\dt \rho  + \nabla\cdot(\rho \bv[\mu])=0.
\en
But, $\bv[\mu]$ being Lipschitz, this solution is unique, and therefore ${\cal F}(\mu_n)$ 
converges to ${\cal F}(\mu)$, which proves the continuity of ${\cal F}$, and ends the proof 
of existence by the Schauder fixed point Theorem.

$\hfill\Box$
\\
We state here without proof some consequences of the previous result:

\begin{cor}\label{5cor}
Let $\rho^0\in {\cal P}(\T^3)$, such that $0< m \leq \rho \leq M$. 
\begin{enumerate}
\item If $\rho^0\in C^{\alpha}, \alpha \in ]0,1]$, 
for $T^*>0$ depending on $\rho^0$,
a solution $\rho(t,x)$ to (\ref{5sg1},\ref{5sg2},\ref{5sg3}) exists 
in $\Linf([0,T^*[, C^{\alpha}(\T^3))$.
\item If $\rho^0 \in W^{1,p}, p>3$, for $T^*>0$ depending on $\rho^0$,
a solution $\rho(t,x)$ to (\ref{5sg1},\ref{5sg2},\ref{5sg3}) exists 
in $\Linf([0,T[, W^{1,p}(\T^3))$.
\item If $\rho^0\in C^{k,\alpha}, \alpha \in ]0,1], k\in \N$, for $T^*>0$ depending on $\rho^0$,
a solution $\rho(t,x)$ to (\ref{5sg1},\ref{5sg2},\ref{5sg3}) exists 
in $\Linf([0,T^*[, 
C^{k,\alpha}(\T^3))$.
\end{enumerate}
Moreover, for these solutions, the velocity field is respectively in $C^{1,\alpha}(\T^3)$, $W^{2,p}(\T^3)$, and $C^{k+1,\alpha}(\T^3)$  on $[0,T^*[$. 
\end{cor}

\section{Uniqueness of solutions to $SG$ with H\"older continuous densities}
\subsection{Result}
Here we prove the following theorem:
\begin{theo}\label{5main4}
Suppose that $\rho^0 \in {\cal P}(\T^3)$ with $0<m \leq \rho^0 \leq M$, 
and belongs to $C^{\alpha}(\T^3)$ for some $\alpha>0$. 
From Theorem \ref{5main3}, for some $T>0$ there exists a solution $\bar \rho$ to $SG$ in $\Linf([0,T],C^\alpha(\T^3))$. 
Then every solution of $SG$ in $\Linf([0,T'],C^{\beta} (\T^3))$ for $T'>0, \beta>0$ with same initial data coincides with $\bar \rho$ on $[0,\inf\{T,T'\}]$.
\end{theo}

{\it Remark 1:} The uniqueness of weak solutions is still an open question.

{\it Remark 2:} Our proof of uniqueness is thus valid in a smaller class of solutions than the one found in the previous section, the reason is the following: during the course of the proof, we will need to solve a Monge-Amp\`ere equation, whose right-hand side is a function of the second derivatives of the solution of another Monge-Amp\`ere equation. 
In Theorem \ref{5wang},  if  $u$ is solution to (\ref{5ma}) with a right hand side satisfying (\ref{5modulus}), although $u\in C^2$,  it 
is not clear that the second derivatives of $u$ satisfy (\ref{5modulus}). Actually, it is even known to be wrong in the case of the Laplacian (for a precise discussion on the subject, the reader may refer to \cite{K}).
However, from Theorem \ref{5caf} below, if $\rho \in C^{\alpha}$ then $u\in C^{2,\alpha}$. 

What we actually need is a continuity condition on the right hand side of (\ref{5ma}) such that the second derivative of the solution $u$
satisfies (\ref{5modulus}). This may be a weaker condition than H\"older continuity, however the proof would
not be affected, therefore it is enough to give it under the present form.

\subsection*{Proof of  Theorem \ref{5main4}} 
Let $\rho_1$ and $\rho_2$ be two solutions of (\ref{5sg1}, \ref{5sg2}, \ref{5sg3}),
in $\Linf([0,T], C^{\beta}(\T^3))$ that coincide at time 0.
Let $X_1, X_2$ be the two corresponding Lagrangian solutions, ({\it i.e.} solutions of
(\ref{5lagrangian},\ref{5lagrangian2})). The velocity field being $C^1$, for all $t\in [0,T]$, $X_1(t,\cdot)$ and $X_2(t,\cdot)$ are both $C^1$ diffeomorphisms of $\Td$.

We call $\bv_1$ (resp. $\bv_2$) the velocity field associated to $X_1$ (resp. $X_2$), $\bv_i(t,x)=[\nabla\Psi_i(t,x)-x]^{\perp}, i=1,2$. 
We have
\be
\dt(X_1-X_2)&=&\bv_1(X_1)-\bv_2(X_2)\\ 
&=&(\bv_1(X_1)-\bv_1(X_2)) + (\bv_1(X_2)-\bv_2(X_2)).
\en
We want to obtain a Gronwall type inequality for $\|X_1-X_2\|_{L^2}$.
Since $\bv_1$ is uniformly Lipschitz in space (from Theorem \ref{5main3}), the first bracket is estimated in $L^2$ norm by $C\|X_1-X_2\|_{L^2}$.
\\
We now need to estimate the second term. We first have that
\be
\int |\bv_1(X_2)-\bv_2(X_2)|^2 = \int \rho_2|\nabla\Psi_1 -\nabla\Psi_2|^2,
\en
and since $\rho_2$ is bounded, we need to estimate
$\|\nabla\Psi_1 -\nabla\Psi_2\|_{L^2}$.  This will be done in the following Proposition:
\begin{prop}\label{5lemme_pour_gronwall}
Let $X_1, X_2$ be mappings from $\Td$ into itself, such that the densities $\rho_i= X_{i\#} dx, i=1,2$ are in $C^{\alpha}(\Td)$ for some $\alpha>0$,  and satisfy $0<m \leq \rho_i \leq M$. Let $\Psi_i, i=1,2$ be convex such that  
\be
\det D^2\Psi_i = \rho_i
\en 
in the sense of Theorem \ref{5defpsi}, {\it i.e.} $\Psi_i=\Psi[\rho_i]$.
Then 
$$\ds \|\nabla\Psi_1-\nabla\Psi_2\|_{L^2} \leq C \|X_1 - X_2\|_{L^2},$$ where $C$ depends on $\alpha$ (the H\"older index of $\rho_i$), $\|\rho_i\|_{C^{\alpha}(\T^d)}$, $m$ and $M$.
\end{prop}

Before giving a proof of this result, we conclude the proof of the Theorem
 \ref{5main4}. 
The Proposition  \ref{5lemme_pour_gronwall} implies immediately that
\be
\|\dt (X_1 -X_2)\|_{L^2} \leq C \|X_1-X_2\|_{L^2},
\en
and we conclude the proof of the Theorem by a standard Gronwall lemma.

$\hfill\Box$ 

\subsection{Energy estimates along Wasserstein geodesics: Proof of Proposition \ref{5lemme_pour_gronwall}.} 

In the proof of this result we will need the following result on optimal transportation of measures by gradient of convex functions:
\begin{theo}[Brenier, \cite{Br1}, McCann, \cite{Mc2}, Cordero-Erausquin, \cite{Co}, Caffarelli,\cite{Ca1}]\label{5caf}
Let $\rho_1$, $\rho_2$ be two probability measures on $\Td$, such that
$\rho_1$ is absolutely continuous with respect to the Lebesgue measure. 
\begin{enumerate}
\item There exists a unique up to a constant convex function $\phi$ such that $\phi - |\cdot|^2/2$ is $\Zd$ periodic, satisfying $\nabla\phi_{\#}\rho_1=\rho_2$.
\item The map $\nabla\phi$ is the solution of the minimization problem 
\beq\label{5opttrans}
\inf_{T_{\#}\rho_1=\rho_2}\int_{\Td}\rho_1(x)|T(x)-x|_{\Td}^2 \ dx,
\enq
and for all $x \in \Rd$,  $|\nabla\phi(x)-x|_{\Td}=|\nabla\phi(x)-x|_{\Rd}$.
\item If $\rho_1$, $\rho_2$ are strictly positive and belong to $C^{\alpha}(\Td)$ for some $\alpha>0$
then $\phi\in C^{2,\alpha}(\Td)$ and satisfies pointwise
\be
\rho_2(\nabla\phi)\det D^2\phi=\rho_1.
\en
\end{enumerate}
\end{theo}

For complete references on the optimal transportation problem (\ref{5opttrans}) and its applications, the reader can refer to \cite{Vi}. 

{\it Remark 1:} the expression $|\cdot|_{\Td}$ denotes the Riemannian distance on the flat torus, whereas $|\cdot|_{\Rd}$ is the Euclidian distance on $\Rd$. The second assertion of point 2 means that, for all $x\in \Rd$, $|\nabla\phi(x)-x| \leq {\rm diam}(\Td)=\sqrt{d}/2$. 

{\it Remark 2:} Here again, note that since $\phi - |\cdot|^2/2$ is periodic, the map $x \mapsto \nabla\phi(x)$ is compatible with the equivalence classes of $\Rd/\Zd$, and therefore is defined without ambiguity on $\Td$.

\subsubsection*{Wasserstein geodesics between probability measures}
In  this part we use results from  \cite{BB1}, \cite{Mc1}.
Using Theorem \ref{5caf}, we consider the unique (up to a constant) convex potential $\phi$ such that
\be
&&\nabla\phi_{\#}\rho_1=\rho_2,\\
&&\phi-|\cdot|^2/2 \ \mbox{is} \  \Zd-\mbox{periodic}.
\en
 We consider, for $\theta\in [1,2]$, $\phi_{\theta}$ defined by
\be
\phi_{\theta}=(2-\theta)\frac{|x|^2}{2} + (\theta-1) \phi.
\en
We also consider, for $\theta\in [1,2]$, $\rho_{\theta}$  defined by
\be
\rho_{\theta}= \nabla\phi_{\theta \  \#}\rho_1.
\en
Then $\rho_{\theta}$ interpolates between $\rho_1$ and $\rho_2$. This interpolation has been introduced in \cite{BB1} and \cite{Mc1} as the time continuous formulation of the Monge-Kantorovitch mass transfer. 
In this construction,  a velocity field $v_{\theta}$ is defined $d\rho_\theta$ a.e. as follows: 
\beq\label{5defvtheta}
\forall f\in C^0(\Td;\Rd), \ \int \rho_\theta v_\theta \cdot f = \int \rho_1f(\nabla\phi_\theta)\cdot \partial_\theta\nabla\phi_\theta.
\enq
It is easily checked that the pair $\rho_\theta, v_\theta$ satisfies
\be
\partial_{\theta}\rho_{\theta} + \nabla\cdot (\rho_{\theta}v_{\theta})=0,
\en 
and for any $\theta\in [1,2]$, we have (see \cite{BB1}):
\be
\demi\int_{\Td} \rho_{\theta}\vert v_{\theta} \vert ^2 = 
 \demi \int_{\Td}\rho_1\vert\nabla\phi(x)-x\vert^2  = W^2_2(\rho_1, \rho_2), 
\en
where $W_2(\rho_1,\rho_2)$ is the Wasserstein distance between $\rho_1$ and $\rho_2$, defined by 
\be
W_2^2(\rho_1,\rho_2) = \inf_{T_{\#}\rho_1=\rho_2}\left\{\int \rho_1(x)\vert T(x)-x\vert_{\Td} ^2\right\}.
\en
The Wasserstein distance can also be formulated as follows:
\be
W^2_2(\rho_1, \rho_2)= \inf_{Y_1, Y_2}\left\{\int_{\Td}|Y_1-Y_2|_{\Td}^2\right\}
\en
where the infimum is performed over all maps $Y_1, Y_2: \Td\mapsto \Td$ 
such that $Y_{i \#}dx=\rho_i, i=1,2$.
From this definition  we have easily
\be
W_2^2(\rho_1,\rho_2) \leq \int \vert X_2(t,a)-X_1(t,a)  \vert ^2 \ da,
\en
and it follows that, for every $\theta\in [1,2]$,
\beq\label{5ineqW2}
\int_{\Td} \rho_{\theta}\vert v_{\theta} \vert ^2 = W^2_2(\rho_1, \rho_2) \leq \|X_2-X_1\|_{L^2}.
\enq

\subsubsection*{Regularity of the interpolant measure $\rho_\theta$}
From Theorem \ref{5caf}, for $\rho_1, \rho_2\in C^{\beta}$ and pinched between the positive postive constants $m$ and $M$ , we know that $\phi \in C^{2,\beta}$ and satisfies
\be
\det D^2\phi = \frac{\rho_1}{\rho_2(\nabla\phi)}.
\en

We now estimate $\rho_{\theta}= \rho_1 [\det D^2\phi_{\theta}]^{-1}$.
From the concavity of $\log(\det (\cdot))$ on symmetric positive matrices, we have
\be
\det D^2\phi_{\theta}&=&\det ((2-\theta)I + (\theta -1)D^2\phi) \\
&\geq& [\det D^2\phi]^{\theta -1}\\ 
&\geq& \frac{m}{M}.
\en
Moreover, since $\phi\in C^2$, $\det D^2\phi_{\theta}$ is bounded by above.
Thus $\rho_{\theta}$ is uniformly bounded away from 0 and infinity, and 
uniformly H\"older continuous. 

\subsubsection*{Final energy estimate}
If we consider, for every $\theta\in[1,2]$, $\Psi_{\theta}$ solution of 
\beq\label{5matheta}
\det D^2\Psi_{\theta}=\rho_{\theta},
\enq
in the sense of Theorem \ref{5defpsiper},
 then
$\Psi_{\theta}$ interpolates between $\Psi_1$ and $\Psi_2$, and $\Psi_{\theta}\in C^{2,\beta}$ uniformly, from the regularity of $\rho_{\theta}$.
We will estimate $\partial_{\theta}\nabla\Psi_\theta$
by differentiating (\ref{5matheta}) with respect to $\theta$:
for $M,N$ two $d\times d$ matrices, $t\in \R$, we recall that 
$$\ds\det(M+tN)=\det M + t \  ({\rm trace} \  M_{co}^t N) + o(t),$$
where $M_{co}$ is the co-matrix (or matrix of cofactors) of $M$.
Moreover, for any $f\in C^2(\Rd;\R)$, if $M$ is the co-matrix of
$D^2f$, it is a common fact that  
\beq\label{5idcom}
\forall j\in [1..d], \sum_{i=1}^d \partial_i M_{ij} \equiv 0.
\enq
Hence, denoting $M_{\theta}$ the co-matrix of $D^2\Psi_{\theta}$,
we obtain that $\partial_\theta\Psi_\theta$ satisfies
\beq\label{5dthetapsi}
\nabla\cdot (M_{\theta}\nabla\partial_{\theta} \Psi_\theta)&=& \partial_{\theta} \rho_{\theta}(t)\nonumber\\
&=& -\nabla\cdot ( \rho_{\theta} v_{\theta} ).
\enq
From the $C^{2,\beta}$ regularity of $\Psi_{\theta}$, $D^2\Psi_{\theta}$ is a $C^{\beta}$ smooth, positive definite matrix, and its co-matrix as well.
Thus the problem (\ref{5dthetapsi}) is uniformly elliptic.
If we multiply by $\partial_\theta\Psi_\theta$, and integrate by parts we obtain
\be
\int \nabla^t\partial_\theta\Psi_\theta\, M_{\theta}\,\nabla\partial_{\theta} \Psi_\theta
 = - \int \nabla\partial_\theta\Psi_\theta \cdot v_\theta \rho_\theta.
\en
Using that $M_\theta\geq \lambda I$ for some $\lambda >0$, and combining with the inequality (\ref{5ineqW2}) above, we obtain
\be
\|\nabla\partial_{\theta}\Psi_{\theta}(t)\|_{L^2}&\leq& \lambda^{-1} \|\rho_{\theta}v_{\theta}\|_{L^2}\\
 &\leq& \lambda^{-1} \|X_2-X_1\|_{L^2}\left(\sup_{\theta}\|\rho_{\theta}\|_{\Linf}\right)^{1/2}.
\en
The constant $\lambda^{-1}$ depends on $m, M, \beta, \{\|\rho_i\|_{C^\beta}, i=1,2\}$, and is thus bounded under our present assumptions.
We have already seen that $\rho_{\theta}$ is uniformly bounded, and 
we finally obtain that
\beq\label{5OKGronwall}
\|\nabla\Psi_1-\nabla\Psi_2\|_{L^2}\leq C \|X_1 -X_2\|_{L^2},
\enq 
this ends the proof of Proposition \ref{5lemme_pour_gronwall}.

$\hfill\Box$.

{\it Remark 1.} In \cite{L2}, the author obtains also (weaker) estimates of the type of Proposition \ref{5lemme_pour_gronwall}, 
for discontinuous densities $\rho_1, \rho_2$.

\section{Uniqueness of solutions to the 2-d Euler equations with bounded vorticity: a new proof}
 This proof adapts easily to the case of 2-d Euler equation with bounded vorticity, giving a new proof of the uniqueness
part in Youdovich's theorem.
We start now from the following system:
\beq
&&\dt \rho + \nabla\psi^\perp\cdot \nabla \rho = 0, \label{5yel1}\\
&&\rho=\Delta\psi,\label{5yel2}\\
&&\rho(t=0)=\rho^0.\label{5yel3}
\enq
For simplicity, we restrict ourselves to the periodic case, {\it i.e.} $x\in \Td$, $\rho, \psi$ periodic, this implies that $\rho$ has total mass equal to 0.  We reprove the following classical result:
\begin{theo}[Youdovich, \cite{Youdo}]\label{5you}
Given an initial data $\rho^0 \in \Linf(\T^2)$ satisfying $\int_{\T^2}\rho^0 = 0$,  there exists a unique solution to (\ref{5yel1}, \ref{5yel2}, \ref{5yel3}) such that  $\rho$ belongs to $\Linf_{loc}(\R^+\times\T^2)$.
\end{theo}
\subsection*{Proof of Theorem \ref{5you}}
We consider two solutions $\rho_1, \psi_1$ and $\rho_2,\psi_2$, such that $\rho_i, i=1,2$ are bounded in $\Linf([0,T]\times \Td)$. 
In this case the velocity fields $\bv_i = \nabla\psi_i^\perp$ both satisfy
\be
\forall (x,y)\in \T^2, \ |x-y| \leq \demi,  \ |\bv_i(x)-\bv_i(y)| \leq C |x-y| \log\frac{1}{|x-y|}.
\en
This implies that  the flows $(t,x)\mapsto X_i(t,x)$ associated to the velocity fields $v_i=\nabla\psi_i^\perp$ are H\"older continuous, and measure preserving. Moreover,  one has, for all $t\in[0,T]$, $\rho_i(t)=X_i(t)_{\#}\rho^0$.

Applying the same technique as before, we need to estimate $\|\nabla\psi_1-\nabla\psi_2\|_{L^2(\T^2)}$ in terms of $\|X_1-X_2\|_{L^2(\T^2)}$. In the present case, the energy estimate of Proposition \ref{5lemme_pour_gronwall} will hold under the weaker assumptions that the two densities are bounded.
\begin{prop}\label{5lpgl}
Let $X_1, X_2$ be  mappings from $\Td$ into itself, let $\rho^0$ be a bounded measure with a density in $\Linf$ with respect to the Lebsgue measure, and with $\int_{\Td} \rho^0 = 0$. Let $\rho_i=X_{i\,\#}\rho^0, i=1,2$.
Let $\psi_i,i=1,2$ be periodic solutions of $\Delta\psi_i = \rho_i, i=1,2$, then we have
\be
\|\nabla\psi_1 - \nabla\psi_2\|_{L^2(\Td)} \leq \left(2\|\rho^0\|_{\Linf}\max \{\|\rho_1\|_{\Linf}, \|\rho_2\|_{\Linf}\} \right)^{1/2} \|X_1-X_2\|_{L^2(\Td)}.
\en
\end{prop}

{\it Remark:} In other words, this proposition shows that for $\rho_1, \rho_2$ bounded, the $H^{-1}$ norm of $\rho_1 -\rho_2$ is controlled by some 'generalized' (since here we have unsigned measures) Wasserstein distance between $\rho_1$ and $\rho_2$. 

To conclude the proof of Theorem \ref{5you}, 
note first that for all $C>0$, we can take $T$ small enough so that $\|X_2-X_1\|_{\Linf([0,T]\times\T^2)} \leq C$. Now we have for the difference $X_1-X_2$, as long as $|X_1-X_2|\leq 1/2$, 
\be
&& \|\dt(X_1-X_2)\|_{L^2} \\
&\leq & \|\nabla\psi_1(X_1)-\nabla\psi_1(X_2)\|_{L^2} +\|\nabla\psi_1(X_2)-\nabla\psi_2(X_2)\|_{L^2}\\  
&\leq & C_1 \| |X_1-X_2| \log(|X_1-X_2|) \|_{L^2} +  C_2 \|X_1-X_2\|_{L^2},
\en  
where, to evaluate the second term of the second line,  we have used the fact that
\be
\|\nabla\psi_1(X_2)-\nabla\psi_2(X_2)\|_{L^2}=\|\nabla\psi_1 - \nabla\psi_2\|_{L^2},
\en
and then applied 
Proposition \ref{5lpgl}.

We just need to evaluate $\||X_1-X_2| \log(|X_1-X_2|) \|_{L^2}$. 
We take $T$ small enough so that $\|X_2-X_1\|_{\Linf([0,T]\times\T^2)} \leq 1/e$
and notice that $x\mapsto x\log^2 x$ is concave for $0\leq x \leq 1/e$, therefore by Jensen's inequality we have
\be
&&\int_{\T^2} |X_2-X_2|^2 \log^2(|X_1-X_2|) \\
&=&\frac{1}{4}\int_{\T^2} |X_2-X_2|^2 \log^2(|X_1-X_2|^2)\\ 
&\leq& \frac{1}{4}\int_{\T^2} |X_2-X_1|^2\, \log^2 \left(\int_{\T^2}|X_2-X_1|^2\right),
\en  
and some elementary computations finally yield  
\be
\dt\|X_2-X_1\|_{L^2} \leq C\|X_2-X_1\|_{L^2} \log\frac{1}{\|X_2-X_1\|_{L^2}}.
\en
The conclusion $X_1 \equiv X_2$ follows then by standard arguments.

\subsection{Energy estimates along Wasserstein geodesic: Proof of Proposition \ref{5lpgl}}
The proof of this proposition is very close to the proof of Proposition \ref{5lemme_pour_gronwall}, and we will only sketch it, insisting on the specific points.
Here the densities $\rho_i$ can not be of constant sign, since their mean value is zero, hence we introduce $\rho^{0,+}$ (resp. $\rho^{0,-}$) the positive (resp. negative) part of $\rho^0$. Then we introduce $\rho_i^\pm = X_{i\#} \rho^{0,\pm}$.
Note that if the mappings $X_i$ were injective, (which is the case in our present situation) we would have $\rho_i^\pm$ that co\"incides with the positive/negative parts of $\rho_i$, but this can be wrong if $X_i$ is not injective. However what remains is that 
$\rho_i=\rho_i^+  - \rho_i^-$.
Now, $\rho_i^{\pm}, i=1,2$ are 4 positive measures of total mass equal to say $M$, with $M<\infty$. 
\subsubsection*{Wasserstein geodesic}
We interpolate between the positive parts $\rho_i^+$, and the negative part is handled in the same way.
As before we introduce the density $\rho_\theta^+(t)$ 
that interpolates between $\rho_1^+(t)$ and $\rho_2^+(t)$.
In this interpolation, we consider $v_\theta^+$ such that 
\be
\partial_\theta \rho_\theta^+ + \nabla\cdot(\rho_\theta^+ v_\theta^+ ) = 0,
\en
and we introduce as well $\rho_\theta^-, v_\theta^-$.
Then $\rho_\theta=\rho_\theta^+ - \rho_\theta^-$ has mean value 0.
Let the potential $\psi_\theta$ be solution to 
\beq\label{5lap}
\Delta \psi_\theta= \rho_\theta.
\enq 
Note that  $\rho_\theta$ has mean value zero therefore this equation is well posed on $\T^2$, moreover $\psi_\theta$ interpolates between $\psi_1$ and $\psi_2$.

\subsubsection*{Bound on the interpolant measure $\rho_\theta$}
Instead of interpolating between two smooth densities, we interpolate between bounded densities, and use the following result from \cite{Mc1}:
\begin{prop}[McCann, \cite{Mc1}]\label{5discvx} 
Let $\rho_\theta^+$ be the Wasserstein geodesic linking $\rho_1^+$ to $\rho_2^+$ defined above. Then, for all $\theta\in [1,2]$,  
\be
\|\rho_{\theta}^+\|_{\Linf}\leq \max \left\{\|\rho_1^+\|_{\Linf}, \|\rho_2^+\|_{\Linf}\right\}.
\en
The same holds for $\rho_i^-,\rho_\theta^-$. 
\end{prop}

{\it Remark:} This property is often referred to as 'displacement convexity'.

\subsubsection*{Energy estimates}
Now by differentiating (\ref{5lap}) with respect to $\theta$, we obtain
\beq\label{5englap}
\Delta\partial_\theta\psi_\theta=\partial_\theta\rho_\theta=-\nabla\cdot(\rho_\theta^+ v_\theta^+ -\rho_\theta^- v_\theta^- ),
\enq
with $v_\theta^{\pm}$ the interpolating velocity defined as in (\ref{5defvtheta}), and satisfying for all $\theta\in [1,2]$,
\be
\int \rho_\theta^{\pm}|v_\theta^{\pm}|^2 =  W^2_2(\rho_1^{\pm}(t), \rho_2^{\pm}(t)).
\en
Multiplying (\ref{5englap}) by $\partial_\theta \psi_\theta$, and integrating over $\theta \in [1,2]$, we obtain
\be
\|\nabla\psi_1-\nabla\psi_2\|_{L^2(\Td)} &\leq& \int_{\theta=1}^2\|\rho_\theta^+ v_\theta^+\|_{L^2}+\|\rho_\theta^- v_\theta^-\|_{L^2}  \\
&\leq&  W_2(\rho_1^+, \rho_2^+)\left(\sup_{\theta}\|\rho_{\theta}^+\|_{\Linf}\right)^{1/2}\\
&+&W_2(\rho_1^-, \rho_2^-)\left(\sup_{\theta}\|\rho_{\theta}^-\|_{\Linf}\right)^{1/2} .
\en
Note that the energy estimate is easier here than in the Monge-Amp\`ere case, since the problem is immediately uniformly elliptic.
 
The mappings $X_i$  satisfy $X_{i\,\#}\rho_0 = \rho_i$, and
$X_{i\, \#}(\rho_0^{\pm})=\rho_i^{\pm}.$ 
Hence, 
\be
W_2^2(\rho_1^{\pm}, \rho_2^{\pm})&\leq& \int \rho_0^{\pm}|X_1-X_2|^2.
\en
Using Proposition \ref{5discvx},  we conclude:
\be
&&\|\nabla\psi_1-\nabla\psi_2\|_{L^2(\Td)} \\
&\leq & 2 \|X_2-X_1\|_{L^2}\Big(\|\rho^0\|_{\Linf} \max\left\{ \|\rho_1\|_{\Linf},\|\rho_2\|_{\Linf} \right\} \Big)^{1/2}.
\en
This ends the proof of Proposition  \ref{5lpgl}. Note that in our specific case, $X_i$ are Lebesgue measure preserving invertible mappings, therefore $\|\rho_i^{\pm}\|_{\Linf}=\|\rho_0^{\pm}\|_{\Linf}$, and the estimate can be simplified in
\be
\|\nabla\psi_1-\nabla\psi_2\|_{L^2(\Td)} \leq 2 \|\rho_0\|_{\Linf}\|X_2-X_1\|_{L^2(\Td)}.
\en

$\hfill\Box$

{\it Remark:} This technique can be used to conclude uniqueness for many  non-linear systems, 
where a transport equation  and an elliptic equation are coupled. The velocity field is the gradient of a
potential satisfying a elliptic equation whose right hand side depends smoothly on the density.
For example, we have uniqueness of solutions to the Vlasov-Poisson system and Euler-Poisson system with bounded 
density and bounded velocity. The Vlasov-Monge-Amp\`ere and Euler-Monge-Amp\`ere systems have also been studied by the author (\cite{BL}, \cite{L1}), and the same technique apply to yield uniqueness for solutions with $C^{\alpha}$ density and bounded velocity. Note however that to enforce uniform ellipticity, we need for the Monge-Amp\`ere equation the density to be bounded by below which is not the case for the Poisson equation.



\section{Convergence to the Euler equation}
\subsection{Scaling of the system}
Here we present a rescaled version of the 2-d $SG$ system and some
formal arguments to motivate the next convergence results.
Here $x\in \T^2, t\in \R^+$ and for $\bv=(v_1,v_2)\in \R^2$, $\bv^{\perp}$ now means $(-v_2,v_1)$.
Introducing $\psi[\rho]=\Psi[\rho] - |x|^2/2$, where $\Psi[\rho]$ is given by Theorem \ref{5defpsiper},
the periodic 2-d $SG$ system now reads 
\be
&&\dt \rho + \nabla \cdot (\rho \nabla\psi^{\perp})=0,\\
&&\det(I+D^2\psi)=\rho.
\en
If $\rho$ is close to one then $\psi$ should be small, and therefore
one may consider the linearization $\det(I+D^2\psi)=1+ \Delta\psi + O(|D^2\psi|^2)$, that 
yields $\Delta \psi \simeq \rho-1$. Thus for small initial data, {\it i.e.} $\rho^0 -1$ small,
one expects $\psi, \mu=\rho-1$ to stay close to a solution of
the Euler incompressible equation  $EI$
\beq
&&\dt \bar\rho + \nabla \cdot (\bar\rho \nabla\bar\phi^{\perp})=0,\label{5euler1}\\
&&\Delta\bar\phi = \bar\rho.\label{5euler2}
\enq
We shall rescale the equation, in order to consider quantities of order one. We introduce the new unknown
\be
&&\rho\ep (t,x) = \frac{1}{\epsilon}(\rho(\frac{t}{\epsilon},x)-1),\\
&&\psi\ep (t,x)= \frac{1}{\epsilon}\psi(\frac{t}{\epsilon},x).
\en
Then we have 
\be
&&\rho( t) = 1 + \epsilon \rho\ep (\epsilon t),\\
&&\Psi[\rho](t)= |x|^2/2 + \epsilon \psi\ep (\epsilon t),
\en
and we define $\phi\ep $ by 
\be
\epsilon\phi\ep  = |x|^2/2 - \Phi[\rho],
\en
so that
\beq\label{5deftphi}
\nabla\phi\ep = \nabla\psi\ep (\nabla\Phi[\rho]).
\enq
Hence, at a point $x\in \T^2$, $\nabla\phi^{\epsilon \perp}$ is the velocity of the associated dual point $\nabla\Phi[\rho](x)$.
The evolution of this quantities is then governed by the system $SG\epu$
\beq
&&\dt \rho\ep  + \nabla\cdot (\rho\ep  \nabla\psi^{\epsilon \perp})=0,\label{5sgep1}\\
&& \det(I+\epsilon D^2\psi\ep )= 1+ \epsilon\rho\ep \label{5sgep2}.
\enq

{\it Remark:} Note that this system admits global weak solutions with initial data any bounded measure $\rho^{\epsilon\,0}$, as long as 
\beq
&&\int_{\T^2} \rho^{\epsilon\,0}=0,\label{5compa1}\\
&&\rho^{\epsilon\,0} \geq -\frac{1}{\epsilon}.\label{5compa2}   
\enq

Note also that if the pair $(\bar\rho,\bar\phi)$ is solution to the $EI$ system (\ref{5euler1}, \ref{5euler2}), so is the pair $\ds\left(\frac{1}{\epsilon}\bar\rho(\frac{t}{\epsilon},x), \  \frac{1}{\epsilon}\bar\phi(\frac{t}{\epsilon},x)\right)$.

We now present the convergence results. We show that solutions of $SG\epu$ converge to solutions of $EI$ in the following sense: if $\rho^{\epsilon\,0}$, the initial data of $SG\epu$, is close (in some sense depending on the type of convergence we wish to show) to a smooth initial data $\bar\rho^0$ for $EI$, then $\rho\ep$ and $\bar\rho$ remain close for some time. This time goes to $\infty$ when $\epsilon$ goes to 0.

We present two different versions of this result: the first one is for weak solutions of $SG\epu$, and the second one is for Lipschitz solutions.

\subsection{Convergence of weak solutions}

\begin{theo}\label{5main5}
Let $(\rho\ep , \psi\ep )$ be a weak solution of the $SG\epu$ system (\ref{5sgep1}, \ref{5sgep2}). Let $(\bar \rho, \bar\phi)$ be a smooth $C^3([0,T]\times\T^2)$ solution of the $EI$ system (\ref{5euler1}, \ref{5euler2}). 
Let $\phi\ep $ be obtained from $\psi\ep $ as in (\ref{5deftphi}), let $H\epu(t)$ be defined by
\be
H\epu(t)=\demi\int_{\T^2}\left|\nabla\phi\ep  -\nabla\bar\phi\right|^2,
\en
then
\be
H\epu(t)\leq (H\epu(0) + C\epsilon^{2/3}(1+t))\exp{Ct}
\en
where $C$ depends on $\sup_{0\leq s\leq t}\{\|D^3\bar\phi(s),D^2\dt\bar\phi(s) \  \|_{\Linf(\T^2)}\}$. 
\end{theo}

{\it Remark 1:} Note that $\nabla\phi^{\epsilon \perp}(t,x)$ is the velocity at point $\nabla\Phi[\rho]=x- \epsilon\nabla\phi\ep $. Thus we compare the $SG\epu$ velocity at point $x- \epsilon\nabla\phi\ep $ (the dual point of $x$) with the $EI$ velocity at point $x$.
Our result allows also to compare the velocities at the same point, by noticing that
\be
G\epu(t)&=&\demi\int_{\T^2}\rho \left|\nabla\psi\ep  -\nabla\bar\phi\right|^2\\
&=&\demi\int_{\T^2}\left|\nabla\phi\ep  -\nabla\bar\phi(x-\epsilon\nabla\phi\ep )\right|^2\\ 
&\leq& C( H\epu(t) + \epsilon^2)
\en
using the smoothness of $\bar\phi$, and if $\bv_{sg\epu},\bv_{ei}$ are the respective velocities of the $SG\epu$ and $EI$ systems, $\ds G\epu=\int_{\T^2}\rho\ep| \bv_{sg\epu}-\bv_{ei}|^2$.

{\it Remark 2:} The expansion $\det(I+D^2\psi)=1+\Delta\psi+O(|D^2\psi|^2)$, used above to justify the convergence relies {\it a priori} on the control of $D^2\psi$ in the sup norm. But in the Theorem \ref{5main5}, the initial data must satisfy 
$\nabla\psi\ep $ close in $L^2$ norm to a smooth divergence free velocity: this condition means that $D^2\psi\ep $ is close in $H^{-1}$ norm to $D^2\bar\phi$, which is smooth. This control does not allow to justify the expansion 
$\det(I+D^2\psi)=1+\Delta\psi+O(|D^2\psi|^2)$, but we see that the result remains valid.
\subsubsection*{Proof of Theorem \ref{5main5}}
In all the proof, we use $C$ to denote any quantity that depends only on $\bar\phi$.
We use the conservation of the energy of the $SG\epu$ system, given by
\be
E(t)=\int_{\T^2}|\nabla\phi\ep |^2. 
\en
This fact, although formally easily justified, is actually not 
so straightforward for weak solutions, and has been proved by F. Otto in an unpublished work. The argument is explained in \cite{BL}. Therefore $E(t)=E_0$.
The energy of the smooth solution of $EI$ is given by 
\be
\int_{\T^2}|\nabla\bar\phi|^2 
\en
and also conserved. 
For all smooth $\theta$, we will use the notation:
\be
<D^2\theta>(t,x)=\int_{s=0}^1 (1-s)D^2 \theta(t,x-s\epsilon \nabla\phi\ep (t,x)).
\en
Thus we have the identity
\beq
\int_{\T^2}\rho\ep \theta &=& \int_{\T^2}\theta(x-\epsilon \nabla\phi\ep )\label{5facile}\\
&=& \int_{\T^2}\theta - \epsilon \int_{\T^2}\nabla\theta\cdot\nabla\phi\ep 
+ \epsilon^2 \int_{\T^2}<D^2\theta>\nabla\phi\ep  \nabla\phi\ep\label{5idtheta} .
\enq
Using the energy bound, the last term is bounded by $\epsilon^2\|D^2\theta\|_{\Linf(\T^2)}E_0$. Then we write
\be
\Dt H\epu(t)&=& -\Dt \int_{\T^2}   \nabla \bar\phi\cdot \nabla\phi\ep.
\en 
Using the identity (\ref{5idtheta}), we have for all smooth $\theta$,
\be
\epsilon \int_{\T^2}\nabla\theta\cdot\nabla\phi\ep=-\int_{\T^2}\rho\ep\theta + \int_{\T^2}\theta + \epsilon^2 \int_{\T^2}<D^2\theta>\nabla\phi\ep  \nabla\phi\ep,
\en
hence, replacing $\theta$ by $\bar \phi$ in this identity, we get 
\be
\Dt H\epu(t)&=& \frac{1}{\epsilon}\Dt \int_{\T^2} [\rho\ep\bar\phi-\bar\phi-
\epsilon^2<D^2\bar\phi> \nabla\phi\ep  \nabla\phi\ep ].
\en
We can suppose without loss of generality that $\ds\int_{\T^2}\bar\phi(t,x) \ dx \equiv 0$. Then if we define
\be
Q\epu(t)=\int_{\T^2} \epsilon<D^2\bar\phi> \nabla\phi\ep  \nabla\phi\ep , 
\en (note that $\vert Q\epu(t)\vert\leq C\epsilon$),  we have
\be
\Dt (H\epu + Q\epu)= \frac{1}{\epsilon}\Dt \int_{\T^2} \rho\ep \bar\phi.
\en
Hence we are left to compute
\be
\frac{1}{\epsilon}\Dt \int_{\T^2} \rho\ep \bar\phi&=&\frac{1}{\epsilon}\int_{\T^2} \dt \rho\ep \bar\phi + \rho\ep \dt \bar\phi\\
&=& \frac{1}{\epsilon}\int_{\T^2}\rho\ep \nabla\psi^{\epsilon\perp} \cdot \nabla \bar\phi -\epsilon \nabla\phi\ep \cdot  \nabla \dt \bar\phi + \epsilon^2<D^2\dt\bar\phi>\nabla\bar\phi \nabla\bar\phi\\
&=& \frac{1}{\epsilon}\int_{\T^2}\rho\ep \nabla\psi^{\epsilon\perp} \cdot \nabla \bar\phi - \int_{\T^2}\nabla\phi\ep \cdot  \nabla \dt \bar\phi + O(\epsilon)\\
&=&T_1+T_2+O(\epsilon),
\en
where at the second line we have used   (\ref{5sgep1}) for the first term and (\ref{5idtheta}) with $\theta = \dt\bar\phi$ for the second and third term. (Remember also that we assume $\int \dt\bar\phi \equiv 0$.)

We will now use the other formulation of the Euler equation: $\bv = \nabla\bar\phi^\perp$ satisfies
\be
\dt\bv+ \bv\cdot\nabla\bv = -\nabla p.
\en
After a rotation of $\pi/2$, this equation becomes: 
\be
\dt\nabla\bar\phi + D^2\bar\phi \nabla \bar\phi^{\perp} = \nabla p^{\perp},
\en 
thus for $T_2$ we have
\be
T_2&=&-\int_{\T^2}\nabla\phi\ep \cdot  \nabla \dt \bar\phi\\ 
&=& \int_{\T^2}\nabla\phi\ep  D^2\bar\phi \nabla\bar\phi^{\perp}.
\en
For $T_1$,  using (\ref{5deftphi}) and (\ref{5idtheta}), we have 
\be
\epsilon T_1&=&\int_{\T^2}\rho\ep \nabla\psi^{\epsilon \perp} \cdot \nabla \bar\phi\\
&=& \int_{\T^2}\nabla\psi^{\epsilon \perp}(x-\epsilon\nabla\phi\ep ) \cdot \nabla \bar\phi(x-\epsilon\nabla\phi\ep )\\
&=& \int_{\T^2}\nabla\phi^{\epsilon \perp}\cdot\nabla\bar\phi  
- \epsilon \nabla\phi^{\epsilon \perp} D^2\bar\phi \nabla \phi\ep   +
\epsilon \Delta
\en
where $\Delta$ is defined by 
\beq\label{5defDelta}
\ds \Delta= \int_{\T^2}\nabla \phi^{\epsilon \perp} \left( D^2\bar\phi -\int_{s=0}^1  D^2\bar\phi(x-s\epsilon \nabla\phi\ep ) \ ds \  \right) \nabla\phi\ep.
\enq
The term $\int_{\T^2}\nabla\phi^{\epsilon \perp}\cdot\nabla\bar\phi$ vanishes identically. Concerning $\Delta$, we claim the following estimate: 
\begin{lemme}\label{5Delta}
Let $\Delta$ be defined by (\ref{5defDelta}), then
\be
|\Delta| \leq C(\epsilon^{\frac{2}{3}}+H\epu),
\en
where $C$ depends on $\|D^3\bar\phi\|_{\Linf}$.
\end{lemme}
 
We postpone the proof of this lemma after the proof of Theorem \ref{5main5}.
We now obtain
\be
\Dt (H\epu(t)+ Q\epu(t))&\leq&\int_{\T^2} (\nabla\bar\phi^{\perp} - \nabla\phi^{\epsilon \perp}) D^2\bar\phi \nabla\phi\ep  + CH\epu +  C\epsilon^{2/3}.
\en
Noticing that for every $\theta:\T^2\mapsto \R$ we have 
\be
\int_{\T^2}\nabla\theta^{\perp} D^2\bar\phi \nabla\bar\phi=\int_{\T^2}\nabla\theta^{\perp} \cdot \nabla (\demi|\nabla\bar\phi|^2)=0,
\en
we find that 
\be
\int_{\T^2} (\nabla\bar\phi^{\perp} - \nabla\phi^{\epsilon \perp}) D^2\bar\phi \nabla\phi\ep  
 = \int_{\T^2}(\nabla\phi^{\perp} - \nabla \bar\phi^{\epsilon \perp}) D^2\bar\phi (\nabla\phi\ep -\nabla\bar\phi),
\en
hence
\be
\Dt (H\epu(t)+ Q\epu(t))&\leq&-\int_{\T^2} (\nabla\phi^{\epsilon \perp} - \nabla \bar\phi^{\perp}) D^2\bar\phi (\nabla\phi\ep -\nabla\bar\phi)+ CH\epu + C\epsilon^{2/3}\\
&\leq& C(H\epu(t) + Q\epu(t) + \epsilon^{2/3})
\en
using that $Q\epu(t) \leq C\epsilon$.
Therefore
\be
H\epu(t)+ Q\epu(t) \leq (H\epu(0)+ Q\epu(0) +C\epsilon^{2/3} t)\exp(Ct)
\en
and finally
\be
H\epu(t)\leq (H\epu(0)+C\epsilon^{2/3} (1+t))\exp(Ct)
\en
and the result follows. Check that the constant $C$ depends only on
\\ 
$\sup_{0\leq s\leq t}\{\|D^3\bar\phi,D^2\dt\bar\phi  \|_{\Linf(\T^2)}\}$.
This ends the proof of Theorem \ref{5main5}

$\hfill \Box$

\subsubsection*{Proof of Lemma \ref{5Delta}} First we show that if $\ds\Theta(R)=\int_{\{\vert\nabla\phi\ep \vert\geq R\}}\vert\nabla\phi\ep \vert^2$, then
\beq\label{5TTA}
\ds\Theta(R)\leq C \int\vert\nabla\phi\ep -\nabla\bar\phi\vert^2 +\frac{C}{R^2}.\enq
Indeed,  $\int \vert\nabla\phi\ep \vert^2\leq C$, implies that
$\ds\mbox{meas}\{\vert\nabla\phi\ep \vert\geq R\}\leq C\frac{1}{R^2}.$
Since $\vert\nabla\bar\phi(t,x)\vert\leq C$ for $(t,x)\in [0,T']\times\Td$, we have
\begin{eqnarray*}
\Theta(R)&\leq&\int_{\{\vert\nabla\phi\ep \vert\geq R\}}\vert\nabla\bar\phi\vert^2+\int_{\{\vert\nabla\phi\ep \vert\geq R\}}\vert\nabla\phi\ep -\nabla\bar\phi\vert^2\\
&\leq& \frac{C}{R^2}+\int\vert\nabla\phi\ep -\nabla\bar\phi\vert^2.
\end{eqnarray*}
Hence (\ref{5TTA}) is proved.

Then, letting
\be
K(x) = D^2\bar\phi -\int_{s=0}^1  D^2\bar\phi(x-s\epsilon \nabla\phi\ep ) \ ds,
\en
we have
\be
&&\Delta\leq C \Theta(R)+ \int_{\vert\nabla\phi\ep \vert\leq R}\vert K(x)|\vert\nabla\phi\ep |^2
\en
with  $\vert K(x)\vert\leq C \epsilon\vert\nabla\phi\ep \vert$ 
thus  
\be
&&\Delta\leq C\epsilon\int_{\vert\nabla\phi\ep \vert\leq R}\vert\nabla\phi\ep \vert^3 + C\Theta(R)\\
&&\leq C\left(\epsilon R \int\vert\nabla\phi\ep \vert^2 +\frac{1}{R^2}+\int \vert\nabla\phi\ep -\nabla\bar\phi\vert^2\right)\\
&&\leq C\left(\epsilon R +\frac{1}{R^2}+\int \vert\nabla\phi\ep -\nabla\bar\phi\vert^2\right)
\en
for all R, so for $R=\epsilon^{-1/3}$ we obtain:
\be
&&\Delta\leq C\epsilon^{2/3}+C\int \vert\nabla\phi\ep -\nabla\bar\phi\vert^2.
\end{eqnarray*}
This proves Lemma \ref{5Delta}

$\hfill \Box$


\subsection{Convergence of strong solutions}
We present here another proof of convergence, that holds for stronger norms.
Let us consider as above the solution $(\bar \rho, \bar \phi)$ to Euler:
\be
&&\dt \bar \rho + \nabla\cdot(\bar\rho\nabla\bar\phi^{\perp}) =0,\\
&&\Delta\bar\phi=\bar\rho,
\en
and we recall the $SG\epu$ system
\be
&&\dt \rho\ep  + \nabla\cdot (\rho\ep  \nabla\psi^{\epsilon \perp})=0,\\
&& \det(I+\epsilon D^2\psi\ep )= 1+ \epsilon\rho\ep .
\en
We have then
\begin{theo}\label{5main6}
Let $(\bar \rho, \bar \phi)$ be a solution of $EI$, such that that $\bar\rho \in C^2_{loc}(\R^+\times\T^2)$.
Let $\rho^{\epsilon\, 0} $ be a sequence of initial data for $SG\epu$ satisfying (\ref{5compa1}, \ref{5compa2}), and such that 
$\displaystyle\frac{\rho^{\epsilon\, 0}-\bar\rho^0}{\epsilon}$ 
is bounded in $W^{1,\infty}(\T^2)$. 
Then there exists a sequence $(\rho\ep , \psi\ep)$ of solutions
to $SG\epu$ that satisfies: for all $T>0$, there exists $\epsilon_T>0$, such that the sequence
\be
\frac{\rho\ep-\bar\rho}{\epsilon}, \frac{\nabla\psi\ep -\nabla\bar\phi}{\epsilon}
\en 
for $0<\epsilon <\epsilon_T$  is uniformly bounded in $\Linf([0,T], W^{1,\infty}(\T^2))$.
\end{theo}

{\it Remark:} In the previous theorem, we obtained estimates in $L^2$ norm, here we obtain estimates in Lipschitz norm.
Estimates of higher derivatives follow in the same way.

\subsubsection*{Proof of Theorem \ref{5main6}}
We expand the solution of $SG\epu$ as the solution of $EI$ plus a small perturbation of order $\epsilon$ and show that this perturbation remains bounded in large norms (at least Lipschitz). 
We first remark the the assumption on $\bar\rho$ implies that
$\forall T>0$, $\bar\phi\in \Linf([0,T]; C^3(\T^2))$.
Let us write
\be
&&\rho\ep  =\bar\rho+ \epsilon\rho_1\\
&&\psi\ep =\bar\phi + \epsilon \psi_1.
\en
Rewritten in terms of $\rho_1, \psi_1$, the $SG\epu$ system reads:
\be
&&\dt\rho_1 + (\nabla\bar\phi + \epsilon\nabla\psi_1)^{\perp}\cdot \nabla\rho_1 = -\nabla\psi_1^{\perp} \cdot\nabla\bar\rho,\\
&&\Delta\psi_1 + \epsilon\textrm{ trace }[ D^2\psi_1 D^2\bar\phi] + \epsilon^2 \det D^2\psi_1 = \rho_1 - \det D^2\bar\phi.
\en
Differentiating the first equation with respect to space, we find the evolution equation for $\nabla \rho_1$:
\beq
&&\dt \nabla \rho_1 + ((\nabla\bar\phi + \epsilon\nabla\psi_1)^{\perp}\cdot \nabla )\nabla \rho_1\nonumber\\
& =& -( D^2\bar\phi  + \epsilon D^2\psi_1) \nabla\rho_1^{\perp} - D^2\psi_1\nabla\bar\rho^{\perp} - D^2\bar\rho \nabla\psi_1^{\perp}\label{5nablarho} .
\enq
We claim that in order to conclude the proof it is enough to have an estimate of the form
\beq\label{5bound}
\|\psi_1(t,\cdot)\|_{C^{1,1}(\T^2)} \leq C(1+\|\rho_1(t,\cdot)\|_{C^{0,1}(\T^2)}),
\enq
where $C$ depends on $\bar\phi$.
Let us admit this bound temporarily, and finish the proof of the theorem:
using (\ref{5bound}) and (\ref{5nablarho}), we obtain 
\be
\Dt \|\nabla \rho_1\|_{\Linf} \leq C(t)(1+ \|\nabla \rho_1\|_{\Linf}+ \epsilon\|\nabla \rho_1\|^2_{\Linf}),
\en
where the constant $C(t)$ depends on the $C^2(\T^2)$ norm of $(\bar\rho(t,\cdot), \bar\phi(t,\cdot))$. This quantity is bounded on every interval $[0,T]$. 

Thus we conclude using Gronwall's lemma  that $\|\nabla\rho_1(t,\cdot)\|_{\Linf(\T^2)}$ remains bounded on $[0,T\epu]$ with $T\epu$ going to $T$ as $\epsilon$ goes to 0.
We then choose $T$ as large as we want, since when $d=2$ the smooth solution to $EI$ is global in time.
From estimate (\ref{5bound}) the $W^{1,\infty}$ bound on $\rho_1$ implies a $W^{2,\infty}$ bound on $\psi_1$.
Then, we remember that 
\be
\rho_1=\frac{\rho\ep-\bar\rho}{\epsilon},  \  \  \ \nabla\psi_1=\frac{\nabla\psi\ep -\nabla\bar\phi}{\epsilon}
\en 
to conclude the proof of Theorem \ref{5main6}.

$\hfill\Box$

\subsubsection*{Proof of the estimate (\ref{5bound})}
We write the equation followed by $\psi_1$ as follows:
\be
\Delta \psi_1= -\textrm{ trace }[\epsilon D^2\psi_1 D^2\bar\phi] - \epsilon^2 \det D^2\psi_1 + \rho_1 - \det D^2\bar\phi.
\en
We recall that
\be
\|f g\|_{C^{2, \alpha}} \leq \|f\|_{C^{2, \alpha}} \|g\|_{C^{2, \alpha}},
\en
hence, using Schauder $C^{2,\alpha}$ estimates for solutions to Laplace equation (see \cite{GT}), we have
\beq\label{5epsibound}
\|\psi_1\|_{C^{2, \alpha}} \leq C_1(1 + \epsilon \|\psi_1\|_{C^{2, \alpha}} + \epsilon^2 \|\psi_1\|_{C^{2, \alpha}}^2),
\enq 
where $C_1$ depends on $\|\bar \phi\|_{C^{2,\alpha}}, \|\rho_1\|_{C^\alpha}$.
The inequality (\ref{5epsibound}) will be satisfied in two cases: either for $\|\psi_1\|_{C^{2, \alpha}}  \leq C_2$ or for $\|\psi_1\|_{C^{2, \alpha}}  \geq C_3\epsilon^{-2}$ where $C_2, C_3$ are positive constants that depend on $C_1$.

Now we show that $\psi\ep$, solution of (\ref{5sgep2}), is bounded in $C^{2,\alpha}$ for $\rho\ep$ bounded in $C^{\alpha}$ norm.
We consider for $t\in [0,1]$ $\psi_t\ep$ the unique up to a constant periodic solution of 
\be
\det(I+\epsilon D^2\psi_t\ep)=1+t\epsilon\rho\ep.
\en
Diiferentiating this equation with respect to $t$, we find
\be
M_{ij}D_{ij}\partial_t \psi_t\ep = \rho\ep,
\en
where $M$ is the co-matrix of $I+\epsilon D^2\psi_t\ep$. From the regularity result of Theorem \ref{5caf}, $M$ is $C^\alpha$ and striclty elliptic. From Schauder estimates, we have
then $\|\partial_t \psi_t\ep\|_{C^{2,\alpha}}\leq C \|\rho\ep\|_{C^{2,\alpha}}$, and integrated over $t\in [0,1]$, we get 
\be
\|\psi\ep\|_{C^{2,\alpha}}\leq C \|\rho\ep\|_{C^{2,\alpha}}.
\en
Hence, since $\psi\ep = \bar\phi+ \epsilon \psi_1$, we have $\psi_1$ bounded by $C/\epsilon$ in $C^{2,\alpha}$. Hence it can not be bigger than $C_3/\epsilon^2$, and  to satsify (\ref{5epsibound}), we must have
\be
\|\psi_1\|_{C^{2, \alpha}}  \leq C_2,
\en
where $C_2$ as above depends on $\|\bar \phi\|_{C^{2,\alpha}}, \|\rho_1\|_{C^\alpha}$. This proves estimate (\ref{5bound}).

$\hfill \Box$

\textbf{Acknowledgment:} The author thanks Mike Cullen for his remarks, and also Yann Brenier, since part
of this work was done under his direction, during the author's PhD thesis. He also thanks Robert McCann and the Fields Institute of Toronto for their hospitality.
\bibliography{biblio}

\def\cprime{$'$} \def\cprime{$'$}
\begin{thebibliography}{10}

\bibitem{BB2}
J.-D. Benamou and Y.~Brenier.
\newblock Weak existence for the semigeostrophic equations formulated as a
  coupled {M}onge-{A}mp\`ere/transport problem.
\newblock {\em SIAM J. Appl. Math.}, 58(5):1450--1461 (electronic), 1998.

\bibitem{BB1}
J.-D. Benamou and Y.~Brenier.
\newblock A computational fluid mechanics solution to the {M}onge-{K}antorovich
  mass transfer problem.
\newblock {\em Numer. Math.}, 84(3):375--393, 2000.

\bibitem{Br1}
Y.~Brenier.
\newblock Polar factorization and monotone rearrangement of vector-valued
  functions.
\newblock {\em Comm. Pure Appl. Math.}, 44(4):375--417, 1991.

\bibitem{Br2}
Y.~Brenier.
\newblock Convergence of the {V}lasov-{P}oisson system to the incompressible
  {E}uler equations.
\newblock {\em Comm. Partial Differential Equations}, 25(3-4):737--754, 2000.

\bibitem{BL}
Y.~Brenier and G.~Loeper.
\newblock A geometric approximation to the {E}uler equations : the
  {V}lasov-{M}onge-{A}mp\`ere equation.
\newblock {\em Geom. Funct. Anal.}, in press.

\bibitem{Ca1}
L.~A. Caffarelli.
\newblock Interior {$W\sp {2,p}$} estimates for solutions of the
  {M}onge-{A}mp\`ere equation.
\newblock {\em Ann. of Math. (2)}, 131(1):135--150, 1990.

\bibitem{Ca0}
L.~A. Caffarelli.
\newblock A localization property of viscosity solutions to the
  {M}onge-{A}mp\`ere equation and their strict convexity.
\newblock {\em Ann. of Math. (2)}, 131(1):129--134, 1990.

\bibitem{Ca3}
L.~A. Caffarelli.
\newblock The regularity of mappings with a convex potential.
\newblock {\em J. Amer. Math. Soc.}, 5(1):99--104, 1992.

\bibitem{Co}
D.~Cordero-Erausquin.
\newblock Sur le transport de mesures p\'eriodiques.
\newblock {\em C. R. Acad. Sci. Paris S\'er. I Math.}, 329(3):199--202, 1999.

\bibitem{CuG}
M.~Cullen and W.~Gangbo.
\newblock A variational approach for the 2-dimensional semi-geostrophic shallow
  water equations.
\newblock {\em Arch. Ration. Mech. Anal.}, 156(3):241--273, 2001.

\bibitem{CuMa}
M.~J.~P. Cullen and H.~Maroofi.
\newblock The fully compressible semi-geostrophic system from meteorology.
\newblock {\em Arch. Ration. Mech. Anal.}, 167(4):309--336, 2003.

\bibitem{CP}
M.~J.~P. Cullen and R.~J. Purser.
\newblock Properties of the {L}agrangian semigeostrophic equations.
\newblock {\em J. Atmospheric Sci.}, 46(17):2684--2697, 1989.

\bibitem{GT}
D.~Gilbarg and N.~S. Trudinger.
\newblock {\em Elliptic partial differential equations of second order}, volume
  224 of {\em Grundlehren der Mathematischen Wissenschaften [Fundamental
  Principles of Mathematical Sciences]}.
\newblock Springer-Verlag, Berlin, second edition, 1983.

\bibitem{K}
J.~Kovats.
\newblock Dini-{C}ampanato spaces and applications to nonlinear elliptic
  equations.
\newblock {\em Electron. J. Differential Equations}, pages No.\ 37, 20 pp.\
  (electronic), 1999.

\bibitem{Li}
J.-L. Lions.
\newblock {\em Quelques m\'ethodes de r\'esolution des probl\`emes aux limites
  non lin\'eaires}.
\newblock Dunod, 1969.

\bibitem{L1}
G.~Loeper.
\newblock The quasi-neutral limit of the {E}uler-{P}oisson and
  {E}uler-{M}onge-{A}mp\`ere systems.
\newblock in preparation.

\bibitem{L2}
G.~Loeper.
\newblock On the regularity of the polar factorisation for time dependent maps.
\newblock {\em Calc. Var. Partial differential equations}, 2004.

\bibitem{LN}
M.~C. Lopes~Filho and H.~J. Nussenzveig~Lopes.
\newblock Existence of a weak solution for the semigeostrophic equation with
  integrable initial data.
\newblock {\em Proc. Roy. Soc. Edinburgh Sect. A}, 132(2):329--339, 2002.

\bibitem{Mc1}
R.~J. McCann.
\newblock A convexity principle for interacting gases.
\newblock {\em Adv. Math.}, 128(1):153--179, 1997.

\bibitem{Mc2}
R.~J. McCann.
\newblock Polar factorization of maps on {R}iemannian manifolds.
\newblock {\em Geom. Funct. Anal.}, 11(3):589--608, 2001.

\bibitem{Vi}
C.~Villani.
\newblock {\em Topics in optimal transportation}, volume~58 of {\em Graduate
  Studies in Mathematics}.
\newblock American Mathematical Society, Providence, RI, 2003.

\bibitem{W}
X.~J. Wang.
\newblock Remarks on the regularity of {M}onge-{A}mpère equations.
\newblock In {\em Proceedings of the International Conference on Nonlinear
  P.D.E.(Hangzhou, 1992)}.

\bibitem{Youdo}
V.~Youdovitch.
\newblock Non-stationary flows of an ideal incompressible.
\newblock {\em Zh. Vych. Mat.}, 3:1032--1066, 1963.

\end{thebibliography}
 
{\begin{flushright}{G. Loeper\\EPFL, SB, IMA\\
 10015 Lausanne\\e-mail: {\sf
gregoire.loeper@epfl.ch}}
\end{flushright}}

\end{document}